\documentclass[11pt,a4paper]{article}
\usepackage{a4wide,amstext, amsmath,amsfonts,amssymb,amsthm,amscd}
\usepackage{stmaryrd}
\usepackage{fancyhdr}
\usepackage{graphicx,color}
\graphicspath{{./Figs/}}

\newcommand{\bbR}{{\mathbb{R}}}

\newcommand{\nablas}{\nabla_\Gamma}
\newcommand{\nablash}{\nabla_{\Gamma_h}}
\newcommand{\Omegaoh}{\Omega_{h,1}}

\newcommand{\bfn}{\boldsymbol n}
\newcommand{\bfbeta}{\boldsymbol \beta}
\newcommand{\bfI}{\boldsymbol I}
\newcommand{\bfP}{\boldsymbol P}
\newcommand{\bfx}{\boldsymbol x}
\newcommand{\bfa}{\boldsymbol a}
\newcommand{\bfb}{\boldsymbol b}

\newcommand{\bfv}{\boldsymbol v}

\newcommand{\mcK}{\mathcal{K}}

\newcommand{\mcF}{\mathcal{F}}
\newcommand{\mcN}{\mathcal{N}}

\newcommand{\mcU}{\mathcal{U}}

\newcommand{\Gammah}{{\Gamma_h}}

\numberwithin{equation}{section}

\newtheorem{rem}{Remark}[section]

\newcommand{\divs}{\text{div}_\Gamma}
\newcommand{\divsh}{\text{div}_{\Gamma_h}}
\newcommand{\vel}{\boldsymbol{\beta}}
\newcommand{\intf}{\Gamma}
\newcommand{\concS}{u_S}

\usepackage{geometry}
\date{}
\begin{document}
\title{{\bf A cut finite element method for coupled bulk-surface problems on time-dependent domains}}
\author{\Large Peter Hansbo$^*$, Mats~G.~Larson$^{\dag}$, Sara Zahedi$^{\ddag}$\\[4mm]
$^*$Department of Mechanical Engineering, J\"onk\"oping University, \\  SE-551~11 J\"onk\"oping, Sweden \\[2mm]
$^{\dag}$Department of Mathematics and Mathematical Statistics, Ume{\aa} University, \\ SE--901~87~~Ume{\aa}, Sweden\\[2mm]
$^{\ddag}$Department of Mathematics, KTH Royal Institute of Technology,\\
SE-100 44 Stockholm, Sweden}
\maketitle

\begin{abstract}
In this contribution we present a new computational method for coupled bulk-surface problems on time-dependent domains. The method is based on a space-time formulation using discontinuous piecewise linear elements in time and continuous piecewise linear elements in space on a fixed background mesh. The domain is represented using a piecewise linear level set function on the background mesh and a cut finite element method is used to discretize the bulk and surface problems. In the cut finite element method the bilinear forms associated with the weak formulation of the problem are directly evaluated on the bulk domain and the surface defined by the level set, essentially using the restrictions of the piecewise linear functions to the computational domain. In addition a stabilization term is added to stabilize convection as well as the resulting algebraic system that is solved in each time step. We show in numerical examples that the resulting method is accurate and stable and results in well conditioned algebraic systems independent of the position of the interface relative to the background mesh.
\end{abstract}


\section{Introduction}
Problems involving phenomena that take place both on surfaces (or interfaces) and in bulk domains occur in a variety of applications in fluid dynamics and biology.  In this paper, we consider a coupled bulk-surface problem modeling the evolution of soluble surfactants. A soluble surfactant is dissolved in the bulk fluid but also exists in adsorbed form on the interface separating two immiscible fluids. Surfactants have a large influence on the dynamics in multiphase flow systems in that they may cause drop-breakup or coalescence due to their ability to reduce the surface tension. They were for example used to lower the surface tension of oil droplets in the 2010 Deepwater Horizon oil spill so that the oil became more soluble in the water. Other examples of applications where the effects of surfactants are important include drug delivery, treatment of lung diseases, and polymer blending~\cite{SSM03}. 

We consider a coupled system of time-dependent convection-diffusion equations describing the concentration of surfactants in the bulk fluid and on the interface. The interface is moving with a given velocity. From a computational point of view, the main challenge is that the differential equations are defined on domains that are evolving with time and that these domains may undergo strong deformations.

A common strategy is to let the mesh conform to the time-dependent domain, see, \emph{e.g.},~\cite{GaTo12, BGN15}. This technique can be made accurate but requires remeshing and interpolation as an interface evolves with time and leads to significant complications when topological changes such as drop-breakup or coalescence occur, especially in three space dimensions. Therefore, computational methods that allow the interface to be arbitrarily located with respect to a fixed background mesh, so called fixed grid methods, have become highly attractive and significant effort has been directed to their development, see, \emph{e.g.},~\cite{CB03,GrRe07,LI06,HaHa02}. In fixed grid methods a strategy for solving the bulk Partial Differential Equation (PDE)  defined on a domain with the interface as boundary is to extend the PDE to the whole computational domain by for example regularized characteristic functions, cf.~\cite{TLLWV09,ChLai14}. Strategies for solving quantities on evolving surfaces are in general developed on the basis of the interface representation technique. In consequence, existing fixed grid methods are usually tightly coupled to the interface representation. Techniques to represent the interface can be roughly divided into two classes: \emph{explicit}\/ representation, \emph{e.g.}, by marker particles~\cite{Pe77}, and \emph{implicit}\/ representation, \emph{e.g.}, by the level set of a higher dimensional function~\cite{OsSet88}. Existing methods using implicit representation techniques generally extend the surfactant concentration to a region embedding the interface, and instead of a surface PDE, a PDE on a higher dimensional domain must be solved for the interfacial surfactant. Several methods have been proposed, based on explicit~\cite{KT11,LTH08,YP98, MT14,KT14, ChLai14} as well as implicit~\cite{AS03,XuZh03,DE10,ESSW11,XYL12} representations. Most work has been done on insoluble surfactants, \emph{i.e.}, surfactants that are only present at the interface without surfactant mass transfer between the interface and the bulk.

In this paper, we present a new computational method for solving coupled bulk-surface problems on time-dependent domains.  The surface PDE is solved on the interface which can be arbitrarily located with respect to the fixed background mesh. The method is accurate and stable and results in well conditioned linear systems independently of how the interface cuts through the background mesh, and the total mass of surfactants is accurately conserved. 

Our strategy is to embed the time-dependent domain where the PDE has to be solved in a fixed background grid, equipped with a standard finite element space, and then take the restriction of the finite element functions to the time-dependent domain. This idea was first proposed for an elliptic problem with a stationary fictitious boundary in~\cite{HaHa02} and for the Laplace--Beltrami operator on a stationary interface in~\cite{ORG09}. It has been extended to other equations with error analysis, for example the Stokes equations involving two immiscible incompressible fluids with different viscosities and with surface tension~\cite{HLZst}, to PDE:s on time-dependent surfaces in, \emph{e.g.},~\cite{HLZch, ORX14, OR14, DER14}, and to stationary coupled bulk-surface problems with linear coupling terms in~\cite{BHLZ,GOR14}.
These types of methods are referred to as cut finite element methods (CutFEM), since the interface cuts through the background grid in an arbitrary fashion. 

We suggest a CutFEM based on a space-time approach with continuous linear elements in space and discontinuous piecewise linear elements in time. The method presented in~\cite{ORX14} is for solving surface PDEs but is also based on a space-time approach with discontinuous elements in time. However, in our approach we add a consistent stabilization term~\cite{Bur10,BHLZadv} which ensures that 1) our method leads to linear systems with bounded condition number, 2) the discretization of  the surface PDE is stable also for convection dominated problems, and 3) the proposed method relies only on spatial discretizations of the geometry at quadrature points in time and results in the same computations as in the case of stationary problems.
In addition, the total mass of surfactants is accurately conserved using a Lagrange multiplier. Numerical results indicate that the method is optimal order accurate (second order); we have also proven optimal order of accuracy for a related stationary coupled bulk-surface problem with a linear coupling term in~\cite{BHLZ}. 

In this paper, we have used the standard level set method to represent the interface, but other interface representation techniques can be used as well.  We have chosen to concentrate on the challenging task of solving the coupled bulk-surface problem on time-dependent domains and will throughout the paper assume that the velocity field is given. 

The remainder of the paper is outlined as follows. In Section~\ref{sec:govequ} we formulate the coupled bulk surface problem.  In Section~\ref{Sec:interf}  we introduce a discrete approximation of the interface and state our assumptions on the geometry. The computational method for the coupled bulk-surface problem modeling soluble surfactants is presented in Section~\ref{sec:method}. In Section~\ref{sec:numexp}, we show numerical examples and in Section~\ref{sec:conc} we summarize our results.

\section{The coupled bulk-surface problem}
\label{sec:govequ} 
\subsection{The domain}
\begin{figure}
\begin{center}
\includegraphics*[width=0.55\textwidth]{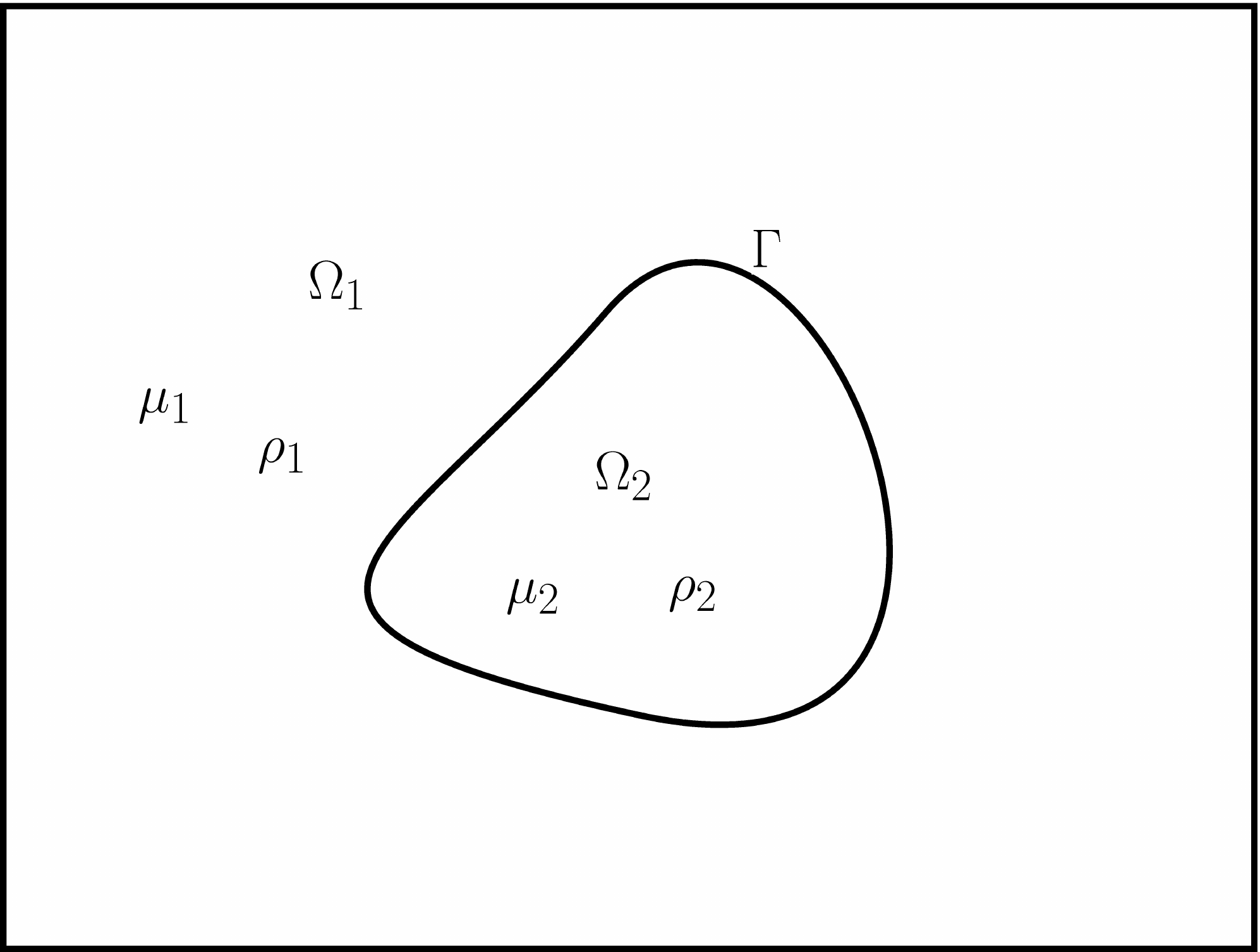} 
\caption{Illustration of the domain $\Omega \in {\bf R}^2$ occupied by two immiscible fluids separated by an interface $\Gamma$.  Immiscible incompressible fluids with density $\rho_i$ and viscosity $\mu_i$ occupy subdomains $\Omega_i(t) \subset \Omega$, $i=1,2$. \label{fig:domain} }
\end{center}
\end{figure}

Let $\Omega$ be an {open bounded domain} in ${\bf R}^d$, $d=2,3$, with convex polygonal boundary $\partial \Omega$ and let $I=[0,T]$ be a time interval. We consider two immiscible incompressible fluids that 
occupy time dependent subdomains $\Omega_i(t) \subset \Omega$, $i=1,2,$ with $t\in I$, such that $\overline{\Omega} = \overline{\Omega}_1(t) \cup \overline{\Omega}_2(t)$ and $\Omega_1(t) \cap \Omega_2(t) = \emptyset$, 
and are separated by a smooth interface defined by $\Gamma(t) = \partial \Omega_1(t) \cap \partial \Omega_2(t)$. See Fig.~\ref{fig:domain} for an illustration in two dimensions. We assume that $\Gamma(t)$ is a simple closed curve (2D) or surface (3D) moving with a given divergence free velocity field $\bfbeta: I \times \Omega \rightarrow {\bf R}^d$ in such 
a way that it does not intersect the boundary $\partial \Omega$ ($\Gamma \cap \partial \Omega = \emptyset$) or itself for any $t\in I$. In applications, the velocity field $\bfbeta$ is typically obtained from the incompressible Stokes or the incompressible Navier--Stokes equations. For simplicity, we further assume that the surfactant is soluble only in the outer fluid phase $\Omega_1(t)$. 

\subsection{The bulk-surface problem}
The model for soluble surfactants is given by a time-dependent convection-diffusion equation on the interface coupled with a time-dependent convection-diffusion equation in the bulk. The concentration of surfactants in the bulk and the concentration of surfactants on the surface are coupled through a nonlinear term which enters as a source term in the PDE for the interfacial surfactant and in a Neumann boundary condition for the bulk concentration. 
More precisely, we consider the following time dependent coupled bulk-surface problem: find 
$u_B: I \times \Omega_1 \rightarrow {\bf R}$ and $u_S: I \times \Gamma \rightarrow {{\bf R}}$ such that 
\begin{alignat}{2}
\partial_t u_B + \bfbeta \cdot \nabla u_B  - \nabla \cdot \left(k_B \nabla u_B \right) &= 0  \qquad &&\text{in $I\times \Omega_1(t)$}   \label{eq:uBP} \\ 
-\bfn\cdot k_B \nabla u_B &= f_{\textrm{coupling}}  \qquad &&\text{on $\Gamma(t)$} \label{eq:BC1} \\ 
-\bfn_{\partial \Omega} \cdot k_B \nabla u_B &= 0 \qquad &&\text{on $\partial \Omega$} \label{eq:BC2} \\
\partial_t u_S + \bfbeta \cdot \nabla u_S +(\divs \bfbeta) u_S  - \divs \left(k_S \nablas u_S\right) &=   f_{\textrm{coupling}}
 &&\text{on $I\times \Gamma(t)$} \label{eq:uSP}
\end{alignat}
with initial condition $u_B(0,\bfx) = u_B^0$ and $u_S(0,\bfx) = u_S^0$ on $\Omega_1(0)$ and $\Gamma(0)$. Here $\partial_t=\frac{\partial}{\partial t}$, $\nabla$ is the usual ${\bf R}^d$ gradient, $\nabla_\Gamma$ is the tangent gradient associated with $\Gamma$ defined by
\begin{equation}
\nabla_\Gamma = \bfP_\Gamma \cdot \nabla, \quad  \bfP_\Gamma = \bfI -\bfn \otimes \bfn
\end{equation}
$\bfn$ is the unit normal vector to $\Gamma$, outward-directed with respect to $\Omega_1$, $\bfn_{\partial \Omega}$ is the outward directed unit normal vector on $\partial \Omega$, $\bfI$\/ is the identity matrix, $\otimes$ denotes outer product $(\bfa \otimes \bfb)_{ij} = a_i b_j$ for any two vectors 
$\bfa$ and $\bfb$), and $k_B$ and $k_S$ are the bulk diffusion and the interfacial diffusion coefficients, respectively. The divergence $\divs \bfv $ on $\Gamma(t)$ of a vector valued function $\bfv$ is defined by
\begin{equation}
\divs \bfv = \text{tr}( \bfv \otimes \nablas ) 
\end{equation}
where $(\bfv \otimes \nablas)_{ij} = (\nablas)_j \bfv_i$.

The exchange of surfactants between the interface and the bulk is modeled with the term $f_{\textrm{coupling}}$. We consider, in particular, the Langmuir model where 
\begin{equation} \label{eq:coupLangm}
f_{\textrm{coupling}}=k_au_B\left(u_S^\infty-u_S \right)-k_du_S
\end{equation} 
Here $u_S^\infty$ is the maximum surfactant concentration on $\Gamma$ and $k_a$ and $k_d$ are adsorption and desorption coefficients, respectively. Examples of other models are 
\begin{alignat}{2}
f_{\textrm{coupling}}&=k_au_B-k_du_S \quad && \textrm{(Henry)} \label{eq:coupHenry} \\
f_{\textrm{coupling}}&=k_au_B\left(1-\frac{u_S}{u_S^\infty} \right)-k_de^{Au_s} u_S \quad  &&(\textrm{Frumkin})\label{eq:coupFrumkin}
\end{alignat}
see for example~\cite{RaFeLi00}. For a numerical study of different isotherms see~\cite{GLS14}.
Using the fact that $\nabla \cdot \bfbeta =0$ we also have the conservation 
law
\begin{equation}\label{eq:conserv}
\int_{\Omega_1(t)} u_B dv +\int_{\Gamma(t)} u_S ds=
\overline{u}_0 \quad \textrm{for $t\geq 0$}
\end{equation}
for the total amount of surfactants on the surface and in the bulk. See the Appendix for the formulation of the transport equations for soluble surfactants in non-dimensional form.

\subsection{Weak form}
First note that we can write $f_{\textrm{coupling}}$, defined in equation \eqref{eq:coupLangm}, in the form 
\begin{equation} \label{eq:coupling}
f_{\textrm{coupling}}={b_B} u_B - {b_S} u_S -b_{BS}u_Bu_S
\end{equation} 
where $b_B=k_au_S^\infty$, $b_S=k_d$, and $b_{BS}=k_a$. We assume that ${b_B}$, ${b_S}$, and $b_{BS}$ are positive constants. For $b_{BS}=0$ we obtain the coupling term in the Henry case, equation \eqref{eq:coupHenry}, with $b_B=k_a$ and $b_S=k_d$.

Let $W=H^1(\Omega_1(t)) \times H^1(\Gamma(t))$. Multiply equation~\eqref{eq:uBP} with a test function $b_Bv_B\in H^1(\Omega_1(t))$ and equation~\eqref{eq:uSP} with a test function $b_Sv_S\in H^1(\Gamma(t))$. After integration by parts and incorporating the boundary conditions~\eqref{eq:BC1} -~\eqref{eq:BC2} we obtain the variational problem: find $u=(u_B,u_S)\in W$ such that
\begin{equation}\label{eq:weakform}
 b_B(\partial_t u_B,v_B)_{\Omega_1(t)} + b_S(\partial_t u_S,v_S)_{\Gamma(t)}+a(t,u,v)-(b_{BS}u_Bu_S, b_Bv_B-b_Sv_S)_{\Gamma(t)} =0 \quad \forall v=(v_B,v_S)\in W
\end{equation}
Here 
\begin{equation}\label{eq:defa}
a(t,u,v)=b_Ba_B(t,u_B,v_B)+b_Sa_S(t,u_S,v_S)+a_{BS}(t,u,v)
\end{equation}
with 
\begin{equation}\label{eq:contforma}
\begin{cases}
a_B(t,u_B,v_B)= (\bfbeta \cdot \nabla u_B, v_B)_{\Omega_1(t)}  + \left(k_B \nabla u_B,\nabla v_B\right)_{\Omega_1(t)} 
\\
a_S(t,u_S,v_S)= (\bfbeta \cdot \nabla u_S,v_S)_{\Gamma(t)} +((\divs \bfbeta) u_S,v_S)_{\Gamma(t)}  + \left(k_S \nablas u_S,\nablas v_S \right)_{\Gamma(t)}
 \\
a_{BS}(t,u,v)=(b_B u_B - b_S u_S,b_B v_B - b_S v_S)_{\Gamma(t)} 
\end{cases}
\end{equation}
Note that $a_{BS}(t,u,v)$ is the linear part of the coupling term $f_{\textrm{coupling}}$.

\section{The level set representation of the interface}\label{Sec:interf}
For $t \in I$, let $\mcU(\Gamma(t))$ be an open neighborhood of $\Gamma(t)$ such that for each $x \in \mcU(\Gamma(t))$ there is a uniquely determined closest point in $\Gamma(t)$. We let 
$\rho(t,\bfx): {\bf R}^d \rightarrow {\bf R}$ be the signed distance function. The exterior unit normal 
$\bfn = \bfn(t,\bfx)$ is the spatial gradient of the signed distance function, $\bfn(t,\bfx ) = \nabla \rho(t,\bfx )$, for $\bfx\in\Gamma(t)$. 

Given a vector field $\bfbeta$ the evolution of the surface $\Gamma(t)$ is governed by the following problem for 
the level set function: find 
$\rho: I \times \Omega \rightarrow {\bf R}$ such that
\begin{equation}\label{eq:levelsetadv}
\rho_t+\bfbeta \cdot \nabla \rho=0, \quad \rho(0,\bfx) = \rho_0
\end{equation}

Let $\mcK_{0,h}$ be a quasiuniform partition of $\Omega$ into shape regular triangles for $d=2$ and tetrahedra for $d=3$ of diameter $h$. We approximate the level set function 
$\rho$ by $\rho_h \in V_{0,h/2}$ where $V_{0,h/2}$ is the space of piecewise linear continuous functions defined on the mesh $\mcK_{0,h/2}$ obtained by refining $\mcK_{0,h}$ uniformly once. $\Gamma_h(t)$ is the zero level set of $\rho_h(t,\bfx)$.  We consider a continuous piecewise linear approximation $\Gamma_h$ of $\Gamma$ such that $\Gamma_h \cap K$ is a linear segment for $d=2$ and is a subset of a hyperplane in $\bbR^3$, for each $K \in \mcK_{0,h/2}$. We assume that for every $t \in I$, $\Gamma_h(t) \subset \mcU(\Gamma(t))$ and that the following approximation assumptions hold:
\begin{equation}\label{eq:geomassumptiona}
\| \rho \|_{L^\infty(\Gamma_h)} \lesssim h^2
\end{equation}
and
\begin{equation}\label{eq:geomassumptionb}
\|\bfn^e - \bfn_h \|_{L^\infty(\Gamma_h)} \lesssim h
\end{equation}
for all $t \in I$. Here $\lesssim$ denotes less or equal up to a positive constant, $\bfn_h$ denotes the piecewise constant exterior unit normal to $\Gamma_h$ with respect to $\Omega_1$ and $\bfn^e$ is the extension of the exact normal to $\Gamma_h$ by the closest point mapping. These assumptions are consistent with the piecewise linear nature of the discrete interface. We define $\Omega_{h,1}$ as the domain enclosed by $\Gamma_h\cup \partial \Omega$. 

Let $0 = t_0 < t_1 <\dots < t_N = T$ be a partition of the time interval $I=[0, T]$ into time steps $I_n = (t_{n-1},t_n]$ of length $k_n = t_n - t_{n-1}$ for $n = 1,2,\dots,N$. To solve the advection equation~\eqref{eq:levelsetadv} we use the Crank--Nicolson scheme in time and piecewise linear continuous finite elements with streamline diffusion stabilization in space. We obtain the method: find $\rho_h^n \in V_{0,h/2}$ 
such that, for $n = 1,2,\dots,N$, 
\begin{align}\label{eq:advdisc}
&(\frac{\rho_h^{n}}{k_n}+\frac{1}{2}\bfbeta^n \cdot \nabla \rho_h^{n},v^n)_\Omega + (\frac{\rho_h^{n}}{k_n}+ \frac{1}{2} \bfbeta^n \cdot \nabla \rho_h^{n},\tau_{SD} \bfbeta^n\cdot \nabla v^n)_\Omega =
\nonumber \\ &\quad =(\frac{\rho_h^{n-1}}{k_n}-\frac{1}{2}\bfbeta^{n-1}\cdot \nabla\rho_h^{n-1}, v^n+\tau_{SD} \bfbeta^n\cdot \nabla v^n)_\Omega  \quad \forall v^n_h \in V_{0,h/2} 
\end{align}
where the streamline diffusion parameter $\tau_{SD}= 2 ( k_n^{-2} + | \bfbeta |^2 h^{-2})^{-1/2}$. 
To keep the level set function a signed distance function, the reinitialization equation, equation (15) in~\cite{SuFa99}, 
can be solved, {\em e.g.}, in the same way as we did with the advection equation in~\eqref{eq:advdisc}.

\section{The space-time cut finite element method}\label{sec:method}
In a space-time formulation the variational formulation is written 
over the space-time domain which is divided into space-time slabs corresponding to each time step.  
We employ piecewise linear trial and test functions that are continuous 
in space but allowed to be discontinuous from one space-time slab to another. For a given time $t_n$ (where $n$ is the time step number) we thus have two distinct solutions, at times $t_n^\pm := \lim_{\epsilon\rightarrow 0}t_n\pm \epsilon$, see, \emph{e.g.}~\cite{Ja78}. Using a suitable weak enforcement of continuity 
at $t_n$ the discrete equations can then be solved one space-time slab at a time.
The discontinuous Galerkin method in time can be compared with implicit finite difference methods and good stability is expected
(the method is related to the first subdiagonal Pad\'e approximation, {\em cf.} Thom\'ee \cite{Th06}).  In space 
we will use a CutFEM, which means that we will use restrictions of standard continuous finite element functions defined on the background grid to our time dependent domain. Since the interface can cut through 
the fixed background grid arbitrarily there is a lack of shape regularity which results in ill-conditioned system matrices. To avoid this problem, we add stabilization terms similar to~\cite{BHLZ}. The same stabilization terms also stabilize the proposed finite element method in case of convection dominated problems.

\subsection{The mesh and finite element spaces}
\begin{figure}\centering
\includegraphics[scale=0.6]{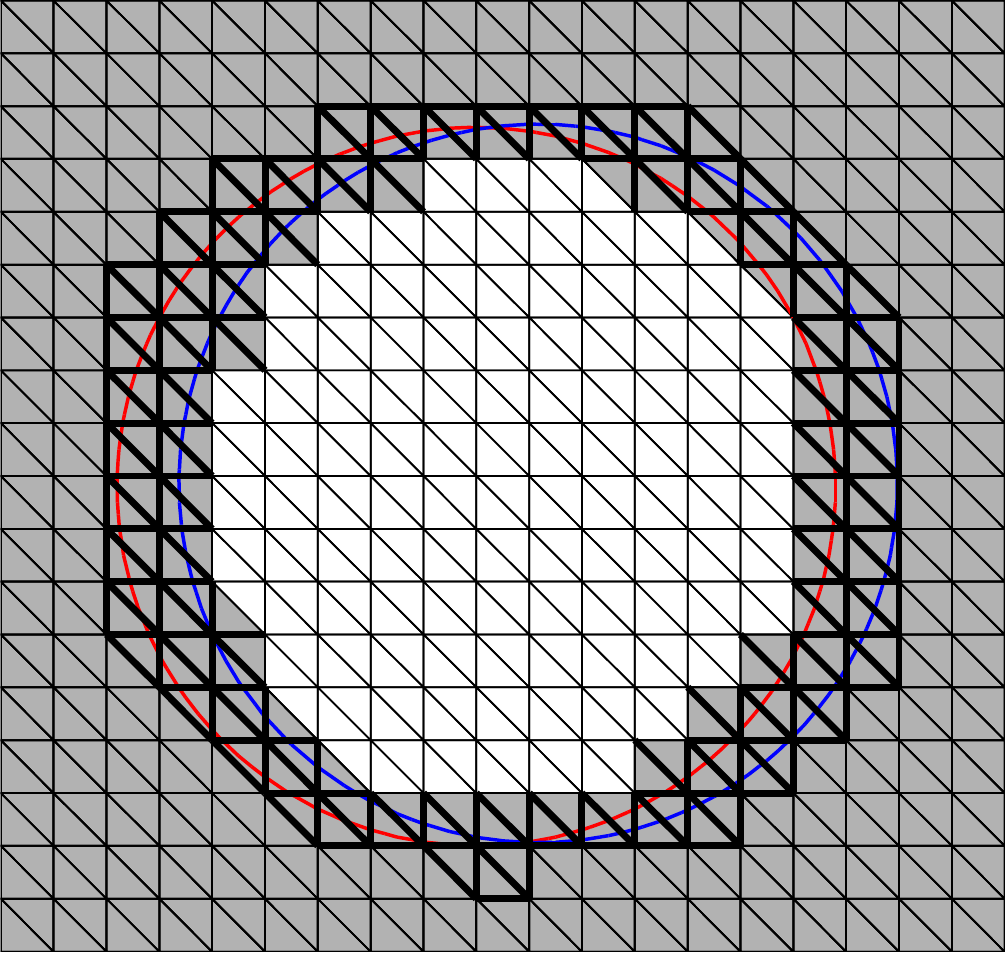} \hspace{2cm}
\includegraphics[scale=0.6]{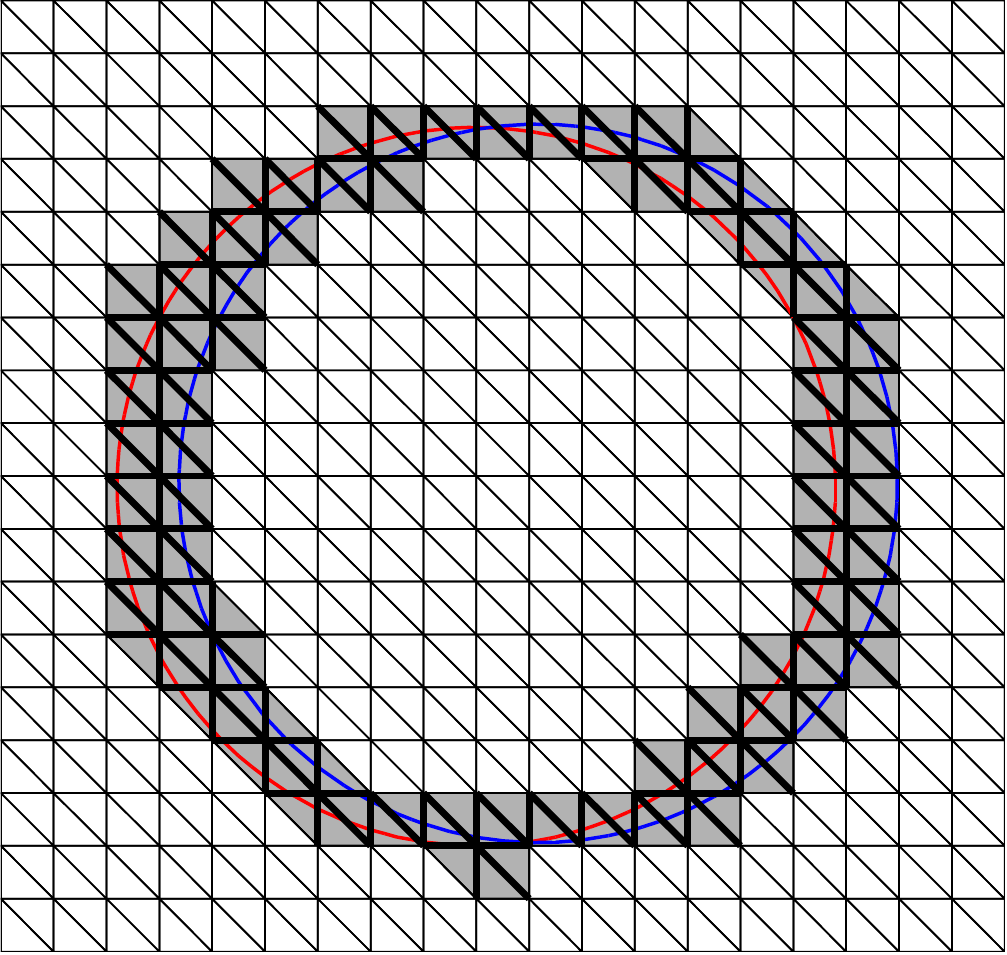}
\caption{Illustration of the sets introduced in Section 4.1. In both figures the blue and the red curves show the position of the interface at the endpoints $t=t_{n-1}$ and $t=t_{n}$ of the time interval $I_n = (t_{n-1},t_n]$.  Left:  the shaded domain shows $\mcN_{B,h}^n$ and edges in $\mcF_{B,h}^n$  are marked with a thick line.  Right: the shaded domain shows  $\mcN_{S,h}^n$ and edges in
$\mcF_{S,h}$ are marked with a thick line.  \label{fig:illust}}.
\end{figure}

We define the following sets of elements
\begin{equation}
\mcK_{B,h}(t) = \{K \in \mcK_{0,h} \,:\, K \cap \Omegaoh(t) \neq \emptyset \},
\qquad 
\mcK_{S,h}(t) = \{K \in \mcK_{0,h} \,:\, K \cap \Gammah(t) \neq \emptyset \}
\end{equation}
and the corresponding sets
\begin{equation}\label{eq:mcnmcs}
\mcN_{B,h}^n = \bigcup_{t \in I_n} \bigcup_{K \in \mcK_{B,h}(t)} K,
\qquad
\mcN_{S,h}^n = \bigcup_{t \in I_n} \bigcup_{K \in \mcK_{S,h}(t)} K  
\end{equation}
For an illustration of the sets in two dimensions see Fig.~\ref{fig:illust}.

Let $P_1(I_n)$ be the space of polynomials of order less or equal to $1$ on $I_n$ and let $V_{0,h}$ be the space of piecewise linear continuous functions defined 
on $\mcK_{0,h}$. On each of the space-time slabs 
$S_{B}^n = I_n \times \mcN_{B,h}^n$ 
and
$S_S^n = I_n \times \mcN_{S,h}^n$ 
we define the spaces
\begin{equation}
V_{B,h}^n = P_1(I_n) \otimes V_{0,h}|_{\mcN_{B,h}^n},
\qquad V_{S,h}^n=P_1(I_n) \otimes V_{0,h}|_{\mcN_{S,h}^n}
\end{equation}
and we let 
\begin{equation}
W_h^n = V_{B,h}^n \times V_{S,h}^n
\end{equation}
Functions in $W_h^n$ take the form
\begin{equation}\label{eq:ansatz}
v(t,\bfx) = (v_B,v_S)=\left(v_{B,0}+v_{B,1}\frac{t-t_{n-1}}{k_n},v_{S,0}+v_{S,1}\frac{t-t_{n-1}}{k_n} \right)
\end{equation}
where $t\in I_n$ and $v_{B,j}$ and $v_{S,j}$, $j=0,1$ can be written as 
\begin{align}\label{eq:funcwh}
v_{B,j}&= \sum_{i=1}^{N_B} \xi_{ij}^B \varphi_i(\bfx)|_{\mcN_{B,h}^n}, \quad 
v_{S,j}= \sum_{i=1}^{N_S}  \xi_{ij}^S \varphi_i(\bfx)|_{\mcN_{S,h}^n} 
\end{align}
Here $\xi_{ij}^B, \xi_{ij}^S \in {\bf R}$ are coefficients, $\varphi_i(\bfx)$ is the standard nodal basis function associated with mesh vertex $i$, $N_B$ and $N_S$ are the number of nodes in $\mcN_{B,h}^n$ and in $\mcN_{S,h}^n$, respectively.

\subsection{The finite element method}\label{sec:weakformST}
Given $u_h(t_{n-1}^-,\bfx)$ and $\overline{u}_0$ (see equation~\eqref{eq:conserv}) find $u_h=(u_B,u_S)\in W_h^n$ and $\lambda \in {\bf R}$, 
such that 
\begin{align}\label{eq:spacetimeform}
A_h^n(u_h,v_h) + J_h^n (u_h,v_h) &+ \lambda\left( (1,v_B)_{\Omega_{h,1}(t_n)} + (1,v_S)_{\Gammah(t_n)} \right) + \mu \left((u_B,1)_{\Omega_{h,1}(t_n)}+ (u_S,1)_{\Gammah(t_n)}\right) \nonumber \\
&=\mu \overline{u}_0,  
\quad \forall v_h \in W_h^n, \mu \in {\bf R}
\end{align}
Here 
\begin{align}
A_h^n (u,v) &= \int_{I_n} b_B(\partial_t u_B,v_B)_{\Omega_{h,1}(t)}\, dt+ \int_{I_n} b_S(\partial_t u_S,v_S)_{\Gammah(t)} \, dt + \int_{I_n} a_h(t,u,v)\, dt 
\nonumber \\
&\quad 
-\int_{I_n}b_{BS} (u_Bu_S, b_Bv_B-b_Sv_S)_{\Gammah(t)} \, dt 
\nonumber \\
&\quad 
+b_B([u_B],v_B(t_{n-1}^+,\bfx))_{\Omega_{h,1}(t_{n-1})} 
+b_S([u_S],v_S(t_{n-1}^+,\bfx))_{\Gammah(t_{n-1})} 
\end{align}
with
\begin{equation}\label{{eq:defa}}
a_h(t,u,v)=b_Ba_{B,h}(t,u_B,v_B)+b_Sa_{S,h}(t,u_S,v_S)+a_{BS,h}(t,u,v)
\end{equation}
and
\begin{equation}\label{eq:contformah}
\begin{cases}
a_{B,h}(t,u_B,v_B)= (\bfbeta \cdot \nabla u_B, v_B)_{\Omega_{h,1}(t)}  + \left(k_B \nabla u_B,\nabla v_B\right)_{\Omega_{h,1}(t)} 
\\
a_{S,h}(t,u_S,v_S)= (\bfbeta \cdot \nabla u_S,v_S)_{\Gammah(t)} +((\divsh \bfbeta) u_S,v_S)_{\Gammah(t)}  + \left(k_S \nablash u_S,\nablash v_S \right)_{\Gammah(t)}
 \\
a_{BS,h}(t,u,v)=(b_B u_B - b_S u_S,b_B v_B - b_S v_S)_{\Gammah(t)} 
\end{cases}
\end{equation}
where $\nablash = \bfP_h \cdot \nabla$ and $\bfP_h = \bfI - \bfn_h \otimes \bfn_h$. Next 
\begin{equation}
J_h^n(u,v)= \int_{I_n} j_h(u,v)\, dt 
\end{equation}
where $j_h(u,v)$ is a stabilizing term of the form
\begin{equation}\label{eq:discformsj}
j_h(v,w)= \tau_B h j_B(v_B,w_B) + \tau_S j_S (v_S,w_S)
\end{equation}
$\tau_B, \tau_S$ are positive parameters and, letting $[x]\vert_F$
denote the jump of $x$ over the face $F$,
\begin{align}
j_B(v_B,w_B) &= \sum_{F \in \mcF_{B,h}} ([\bfn_F\cdot \nabla v_B],[\bfn_F\cdot \nabla w_B])_F
\\
j_S(v_S,w_S) &= \sum_{F \in \mcF_{S,h}} ([\bfn_F \cdot \nabla v_S],[\bfn_F\cdot \nabla w_S])_F
\end{align}
with $\mcF_{S,h}$ the set of internal faces (i.e. faces with two 
neighbors) in $\mcN_{S,h}^n$ and $\mcF_{B,h}$ the set of faces that 
are internal in  $\mcN_{B,h}^n$ and also belong to an element in $\mcN_{S,h}^n$, see Fig.~\ref{fig:illust}. 
Finally, note that we use a Lagrange multiplier to impose the condition~\eqref{eq:conserv}.

\begin{rem}
The stabilization terms $j_B$ and $j_S$ that appear in the method are consistent and are needed to control the condition number of the system matrix so that the resulting algebraic system is well conditioned independently of the position of the interface relative to the mesh. The stabilization term $j_S$ also ensures a stable discretization of the surface PDE in case the problem is convection dominated~\cite{BHLZadv}. However, one may need to apply the stabilization term $j_B$ on all faces that are internal in  $\mcN_{B,h}^n$ to guarantee a stable discretization of the bulk PDE in case the problem is convection dominated. One may also apply the bulk stabilization only on  faces in  $\mcF_{B,h}$ and add other stabilization methods like for example SUPG~\cite{BrHu82} when the problem is convection dominated. 
\end{rem}
  
\begin{rem}
The proposed CutFEM method with continuous linear elements in 
space and discontinuous piecewise linear elements in time is 
second order accurate both in space and time. 
\end{rem}
\begin{rem} \label{rem:dghighorder} We may also consider higher 
order discontinuous Galerkin method in time based on polynomials 
of order $p$.  The optimal order of convergence for a parabolic problem 
is  $p+1$ inside the intervals $I_n$  and $2p+1$ in the nodes $t_n$, 
see \cite{Th06} for further details. Note, however, that in order to achieve higher order convergence in time both the transport 
equation \eqref{eq:levelsetadv} for the distance function $\rho$  and 
the bulk-surface problem (\ref{eq:uBP}-\ref{eq:uSP}) for the concentrations $u_B$ and $u_S$ must be discretized  with a higher order method.
\end{rem}

\begin{rem} 
We may consider the same method without the Lagrange multipliers $\lambda$ and $\mu$. This method is also 
of optimal convergence order but the conservation of the total 
mass of surfactants is lost. See also Example 1 in Section~\ref{sec:numexp}. Strong imposition of the conservation 
law using Lagrange multipliers essentially compensates for 
numerical errors, such as the error in the area of the surface 
and in the volume of the bulk domain, during each time step.
\end{rem}

\subsection{Implementation}
\subsubsection{Newton's method}
Since the bulk and surface surfactant forms are coupled 
through a nonlinear term, see \eqref{eq:coupLangm}, the proposed method~\eqref{eq:spacetimeform} leads 
to a nonlinear system of equations in each time step, 
which we solve using Newton's method. To formulate 
Newton's method we define the residual $F$ and the 
Jacobian $DF$ as follows 
\begin{align}\label{eq:F}
F(u,\lambda)&= \int_{I_n}b_B (\partial_t u_B,v_B)_{\Omega_{h,1}(t)}\, dt+ \int_{I_n} b_S(\partial_t u_S,v_S)_{\Gammah(t)} \, dt + \int_{I_n} a_h(t,u,v)\, dt 
\nonumber \\
&\quad 
-\int_{I_n}b_{BS} (u_Bu_S, b_Bv_B-b_Sv_S)_{\Gammah(t)} \, dt 
\nonumber \\
&\quad 
+b_B([u_B],v_B(t_{n-1}^+,\bfx))_{\Omega_{h,1}(t_{n-1})} 
+b_S([u_S],v_S(t_{n-1}^+,\bfx))_{\Gammah(t_{n-1})} 
+\int_{I_n}j_h(u,v)\, dt 
\nonumber \\
&\quad 
+\lambda\left( (1,v_B)_{\Omega_{h,1}(t_n)} + (1,v_S)_{\Gammah(t_n)} \right) + \mu \left((u_B,1)_{\Omega_{h,1}(t_n)}+ (u_S,1)_{\Gammah(t_n)}\right) -\mu \overline{u}_0
\end{align}
\begin{align}\label{eq:DF}
DF(u,\lambda)(w,\hat{\lambda}) &= \int_{I_n} b_B(\partial_t w_B,v_B)_{\Omega_{h,1}(t)}\, dt+ \int_{I_n} b_S(\partial_t w_S,v_S)_{\Gammah(t)} \, dt + \int_{I_n} a_h(t,w,v)\, dt 
\nonumber \\
&\quad 
-\int_{I_n}b_{BS} (w_Bu_S, b_Bv_B-b_Sv_S)_{\Gammah(t)} \, dt -\int_{I_n}b_{BS} (u_Bw_S, b_Bv_B-b_Sv_S)_{\Gammah(t)} \, dt 
\nonumber \\
&\quad 
+b_B(w_B,v_B(t_{n-1}^+,\bfx))_{\Omega_{h,1}(t_{n-1})} 
+b_S(w_S,v_S(t_{n-1}^+,\bfx))_{\Gammah(t_{n-1})} 
+\int_{I_n}j_h(w,v)\, dt 
\nonumber \\
&\quad 
+\hat{\lambda}\left( (1,v_B)_{\Omega_{h,1}(t_n)} + (1,v_S)_{\Gammah(t_n)} \right) + \mu \left((w_B,1)_{\Omega_{h,1}(t_n)}+ (w_S,1)_{\Gammah(t_n)}\right)
\end{align}
With this notation the nonlinear problem resulting from (\ref{eq:spacetimeform}) takes the form: 
find $u_h\in W_h^n$ and $\lambda \in {\bf R}$ such that 
$F(u_h,\lambda)=0$, and the corresponding Newton's method reads: 
\begin{enumerate}
\item[1.] Choose initial guesses $u_{h,0}$ and $\lambda_0$
\item[2.]  while  $||(w,\hat{\lambda})||>$ tol
\begin{itemize}
\item Solve: $DF(u_{h,0},\lambda_0)(w,\hat{\lambda}) =F(u_{h,0},\lambda_0)$
\item Update $u_{h,0}$: $u_{h,0}=u_{h,0}-w$ and $\lambda_0$: $\lambda_0=\lambda_0-\hat{\lambda}$
\end{itemize}
\end{enumerate}
For $t\in I_n$ we choose the initial guess $u_{h,0}$ to be the solution at $t_{n-1}^- $ i.e. $u_{h,0}(t,\bfx)=u_h(t_{n-1}^-,\bfx)$.

\subsubsection{Assembly of the bilinear forms using quadrature in time}
Recall from equation \eqref{eq:ansatz} that a function $v \in W_h^n$ can be written as
\begin{equation}
v(t,\bfx) = (v_B,v_S)=\left(v_{B,0}+v_{B,1}\frac{t-t_{n-1}}{k_n},v_{S,0}+v_{S,1}\frac{t-t_{n-1}}{k_n} \right)
\end{equation}
where $t\in I_n$ and $v_{B,0}$ and $v_{B,1}$ are functions in $V_{0,h}|_{\mcN_{B,h}^n}$  (the space of restrictions to $\mcN_{B,h}^n$ of functions in $V_{0,h}$) and  $v_{S,0}$ and $v_{S,1}$ are functions in $V_{0,h}|_{\mcN_{S,h}^n}$  (the space of restrictions to $\mcN_{S,h}^n$ of functions in $V_{0,h}$). Thus, for example the first term in $DF(u,\lambda)(w,\hat{\lambda})$ can be written as
\begin{equation}
\int_{I_n} b_B(\partial_t w_B,v_B)_{\Omega_{h,1}(t)}\, dt=
\int_{I_n} \frac{b_B}{k_n}(w_{B,1},v_{B,0})_{\Omega_{h,1}(t)}\, dt+
 \int_{I_n} b_B\frac{t-t_{n-1}}{k_n^2}(w_{B,1},v_{B,1})_{\Omega_{h,1}(t)} \, dt
\end{equation}
Our approach is to first use a quadrature rule in time and for each 
of the quadrature points compute the integrals in space. Since 
the geometry changes in time the contributions to the Jacobian 
and the residual are computed in each quadrature point for the 
current geometry at that point in time. Using a quadrature formula 
in time with quadrature weights $\omega_q^n$ and 
quadrature points $t_q^n$, $q=1,\dots,n_q$, where $n_q$ is the 
number of quadrature points, the first term in $DF(u,\lambda)(w,\hat{\lambda})$ is approximated by
\begin{equation}
\int_{I_n} b_B(\partial_t w_B,v_B)_{\Omega_{h,1}(t)}\, dt\approx
\sum_{q=1}^{n_q} \omega_q^n\frac{b_B}{k_n}(w_{B,1},v_{B,0})_{\Omega_{h,1}(t_q^n)} +
\sum_{q=1}^{n_q} \omega_q^n b_B\frac{t_q^n-t_{n-1}}{k_n^2}(w_{B,1},v_{B,1})_{\Omega_{h,1}(t_q)}
\end{equation}
The other terms in $F(u,\lambda)$ and $DF(u,\lambda)(w,\hat{\lambda})$ 
are treated in the same way. In the numerical examples in Section~\ref{sec:numexp} we use both the trapezoidal rule and Simpson's rule for the time integration.
Choosing the trapezoidal rule for the time integration corresponds to choosing $n_q=2$, quadrature points $t_1^n=t_{n-1}$ and $t_2^n=t_{n}$, and quadrature weights $\omega_1^n=\omega_2^n=\frac{k_n}{2}$. Choosing Simpson's quadrature rule corresponds to choosing 
$n_q=3$, quadrature points $t_1^n=t_{n-1}$, $t_2^n=\frac{t_{n-1}+t_n}{2}$, and $t_3^n=t_{n}$, and quadrature weights $\omega_1^n=\omega_3^n=\frac{k_n}{6}$ and $\omega_2^n=\frac{4k_n}{6}$. For higher 
order discontinuous Galerkin methods in time, see Remark \ref{rem:dghighorder}, it is convenient to 
use the Gauss-Lobatto quadrature rule, which includes the 
endpoints of the interval as quadrature points and is exact for 
polynomials of order $2n_q-3$ and thus have an error term of 
the form $O(k_n^{2n_q-2})$. The optimal local order in a discontinuous 
Galerkin method in time is $2p+2$, where $p$ is the order of polynomials, and therefore it is natural to choose $n_q = p + 2$. When using a quadrature rule that includes the endpoints 
some computations can be reused when passing from one 
space-time slab to another.

For each of the quadrature points $t_q^n$ in the time interval $I_n$  
we compute the discrete surface $\Gamma_h(t_q^n)$ defined as the 
zero level set of the approximate distance function $\rho_h(t_q^n)$. The intersection $\Gamma_h(t_q^n) \cap K$ is then planar, since $\rho_h(t_q^n)$ is piecewise linear, and we can therefore easily compute the contribution to the stiffness matrix using a quadrature rule for a line segment in two dimensions and triangles in three dimensions.  
The contribution from integration on $\Omega_{h,1}(t_q^n) \cap K$ is divided into contributions on one or several triangles in two dimensions and tetrahedra in three dimensions depending on how the interface cuts element $K$.

Note that when we use Simpson's quadrature rule we need to compute $\rho_h(t_n)$, $\rho_h(t_{n+1})$, and $\rho_h(t_{n+1/2})$ with $t_{n+1/2}=\frac{t_{n-1}+t_n}{2}$. Therefore we use half the time step size i.e. $k/2$ when evolving the interface (solving the advection equation for the level set function) while the coupled bulk-surface problem is solved with time step size $k$.

Finally, we obtain a $2(N_B+N_S)+1 \times 2(N_B+N_S)+1$ 
linear system of equations 
\begin{equation}
DF(u_{h,0},\lambda_0)(w,\hat{\lambda}) =F(u_{h,0},\lambda_0)
\end{equation}
for $\hat{\lambda} \in {\bf R}$ and
\begin{equation}
w=\left(\begin{array}{l}
w_{B,0} \nonumber \\
w_{S,0}\nonumber \\
w_{B,1} \nonumber \\
w_{S,1} \nonumber \\
\end{array} \right)
\end{equation}
We use a direct solver to solve this linear system of equations. 


\subsubsection{Different models of the bulk-surface coupling} 
In the proposed method it is straightforward to account for other forms of the coupling term $f_{\textrm{coupling}}$. We now briefly explain how the proposed method should be modified when other examples for $f_{\textrm{coupling}}$ are used. 

The variational formulation and thus $F(u_h,\lambda)$ contains the term $\int_{I_n} (f_{\textrm{coupling}}, b_Bv_B-b_Sv_S)_{\Gammah(t)}\, dt$. The linear part of $f_{\textrm{coupling}}$ is contained in $a_{BS,h}$ and the fourth term in $F(u_h,\lambda)$ (equation~\eqref{eq:F}) is the nonlinear part of $f_{\textrm{coupling}}$. Thus, these two terms and consequently the Jacobian $DF$ changes if the coupling term changes.

For example, in the Henry case~\eqref{eq:coupHenry} $a_{BS,h}$ is as before but the fourth term in $F(u_h,\lambda)$ vanishes since $b_{BS}=0$ and the problem is linear. In the Frumkin case $f_{\textrm{coupling}}$ is of the form 
\begin{equation}
f_{\textrm{coupling}}=b_Bu_B-b_Se^{Au_s}u_S-b_{BS}u_Bu_S
\end{equation}
 Hence $a_{BS,h}=(b_Bu_B, b_Bv_B-b_Sv_S)_{\Gammah(t)}$ and the fourth term in equation~\eqref{eq:F} is replaced by  
\begin{equation}
-\int_{I_n} (b_Se^{Au_s}u_S+b_{BS}u_Bu_S, b_Bv_B-b_Sv_S)_{\Gammah(t)}\, dt
\end{equation}
and the term 
\begin{equation}
-\int_{I_n} (b_S(Ae^{Au_S}u_S+e^{Au_S})w_S, b_Bv_B-b_Sv_S)_{\Gammah(t)}\, dt
\end{equation}
has to be added to equation~\eqref{eq:DF}.

\begin{rem}Another approach to assembly of the discrete problem, 
used for example in~\cite{CL15}, is to explicitly construct three dimensional triangulations, in case of two space dimensions, or four dimensional triangulations, in case of three space dimensions, of the space-time subdomains. However, in comparison our approach is 
very simple to implement since it only relies on spatial discretizations 
of the geometry at the quadrature points in time, which is the same computation as in the case of a stationary problem. 
\end{rem}

\begin{rem} The implementation of the method is straightforward also 
in three spatial dimensions, since the use of quadrature in time essentially reduces the geometric computations to the corresponding stationary  problem. A detailed study of a linear coupled bulk-surface problem, including an implementation in 3D, estimates of the error, and estimates of the condition number, is presented in \cite{BHLZ}. Due to 
the stabilization terms the condition number is always under control 
enabling the use of efficient linear algebra solvers. 
\end{rem}

\begin{rem}\label{rem:nonzerorhs}
If a coupled bulk-surface problem is solved where the bulk and surface PDEs, equation~\eqref{eq:uBP} and~\eqref{eq:uSP}, have nonzero right handsides $f_B(t,\bfx)$ and $f_S(t,\bfx)$ one has to add the terms $\int_{I_n} b_B(f_B,v_B)_{\Omega_{h,1}(t)}\, dt$ and $\int_{I_n} b_S(f_S,v_S)_{\Gammah(t)}\, dt$ to $F(u,\lambda)$. This is the case in Example 1 in Section~\ref{sec:numexp}.
\end{rem}

\section{Numerical examples}\label{sec:numexp}
In all the computations in this section we use a uniform underlying mesh $\mcK_{0,h}$ consisting of triangles of size $h$ and a constant time step size of the form $k = C h$. The stabilization constants $\tau_B$ and $\tau_S$ in the stabilization term $j_h$ are $10^{-2}$.

We consider examples from~\cite{XuZh03, KT14, ChLai14} and we also formulate one example for which we know the exact solution. In the examples from~\cite{KT14, ChLai14} the coupled bulk-surface problem is formulated in non-dimensional form. This leads to some minor changes in the weak formulation, see the Appendix.

In this section we show that the proposed method is second order accurate. We study the convergence at time $t=0.5$. In case the exact solution is known we measure the order of convergence by studying $\| (u_{B,\textrm{exact}} - u_{B,h}) \|_{\Omega_{h,1}(0.5)}$  and $\| (u_{S,\textrm{exact}} - u_{S,2h}) \|_{\Gammah(0.5)}$ for different mesh sizes $h$. When the exact solution is not known we measure the order of convergence by using consecutive refinements of the underlying mesh and study $\| (u_{B,h} - u_{B,2h}) \|_{\Omega_{h,1}(0.5)}$  and $\| (u_{S,h} - u_{S,2h}) \|_{\Gammah(0.5)}$. This is also how the convergence is studied in~\cite{KT14, ChLai14}. 

Note that in the computation of $\| (u_{S,h} - u_{S,2h}) \|_{\Gammah(0.5)}$ the concentration $u_{S,2h}(0.5)$ needs to be evaluated at quadrature points on $\Gammah(0.5)$ which do not have to lie on $\Gamma_{2h}(0.5)$. Thus, we need to extend the concentration $u_{S,2h}(0.5)$ out from $\Gamma_{2h}(0.5)$ to $\Gammah(0.5)$. There are different ways of doing this extension. 
However, analogously to our analysis in~\cite{BHLZ}, for each quadrature point on $\Gammah(0.5)$ (used in the computation of $\| (u_{S,h} - u_{S,2h}) \|_{\Gammah(0.5)}$) we find the closest point on $\Gamma_{2h}(0.5)$ and evaluate $u_{S,2h}$ at that point. We use the same approach in the computation of $\| (u_{B,h} - u_{B,2h}) \|_{\Omega_{h,1}(0.5)}$.  

We will also show in this section that the condition number of the algebraic system of equations is bounded independently of how the interface cuts the underlying mesh and that the total mass of surfactants can be conserved accurately. 

\subsection{Example 1}
To study the convergence of our numerical method we now consider a coupled bulk-surface problem where the exact solution is given. The interface is initially a circle centered at $(0.5,0.22)$ with radius $r_0=0.17$ and the velocity field $\bfbeta=(\pi(0.5-y),\pi(x-0.5))$. The interface moving with this velocity field is at time $t$ a circle with radius $r_0=0.17$ centered at 
\begin{align}\label{eq:xcyc}
x_c&=0.5+0.28\sin(\pi t) \nonumber \\
y_c&=0.5-0.28\cos(\pi t) 
\end{align}
We choose the computational domain to be $[0, \ 1]\times [0, \ 1]$. The bulk diffusion and the interfacial diffusion coefficients are $k_B=0.01$ and $k_S=1$.
The exact solution is
\begin{align}
u_B&=0.5+0.4\cos(\pi x)\cos(\pi y)\cos(2\pi t) \nonumber \\
u_S&=\frac{u_B +
\frac{\pi}{250}\sin(\pi x)\cos(\pi y)\cos(2\pi t)n_1+
\frac{\pi}{250}\cos(\pi x)\sin(\pi y)\cos(2\pi t)n_2}
{1.5+0.4\cos(\pi x)\cos(\pi y)\cos(2\pi t)}
\end{align} 
where  
\begin{align}
n_1&=\frac{(x-x_c)}{\sqrt{(y-y_c)^2+(x-x_c)^2}} \nonumber \\
n_2&=\frac{(y-y_c)}{\sqrt{(y-y_c)^2+(x-x_c)^2}}
\end{align}
and $x_c$ and $y_c$ are given in equation~\eqref{eq:xcyc}. The function $u_B$  satisfies the interface and boundary conditions \eqref{eq:BC1} and \eqref{eq:BC2} but the bulk and surface PDEs~\eqref{eq:uBP} and~\eqref{eq:uSP} are satisfied with right hand sides $f_B$ and $f_S$, respectively. For the implementation, see Remark \ref{rem:nonzerorhs}. 

\begin{figure}\centering
\includegraphics[scale=0.395]{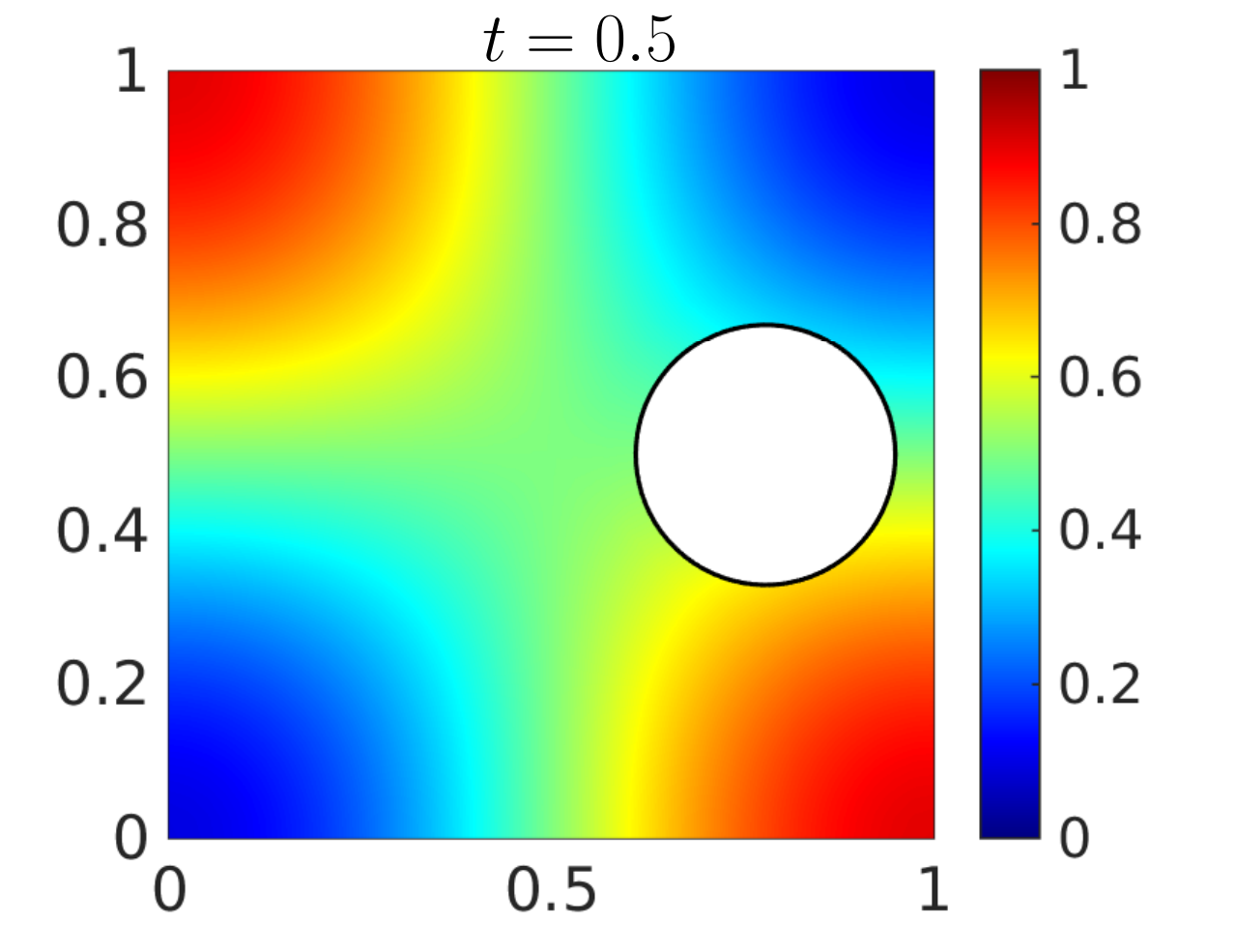} \hspace{0.1cm}
\includegraphics[scale=0.395]{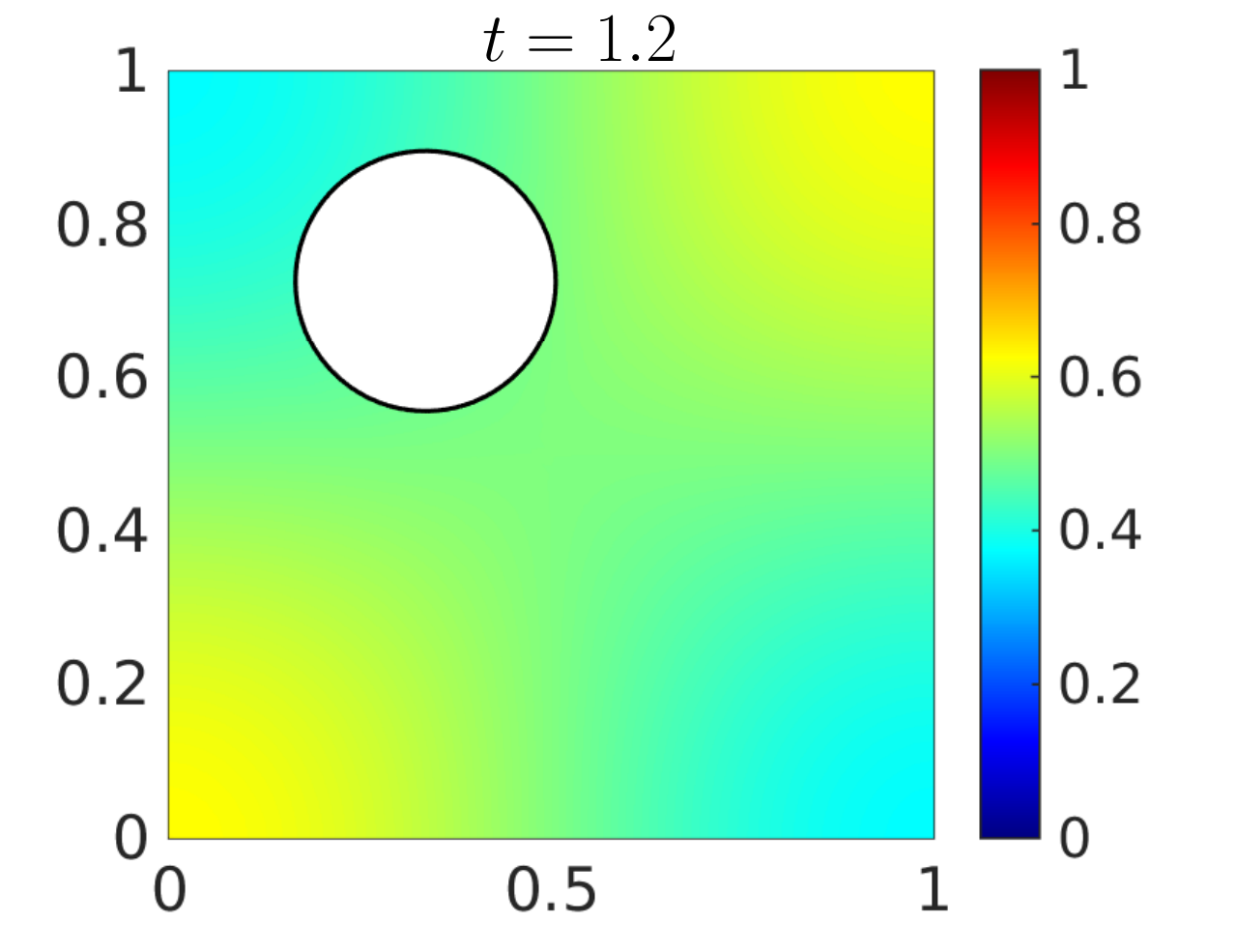}
\includegraphics[scale=0.395]{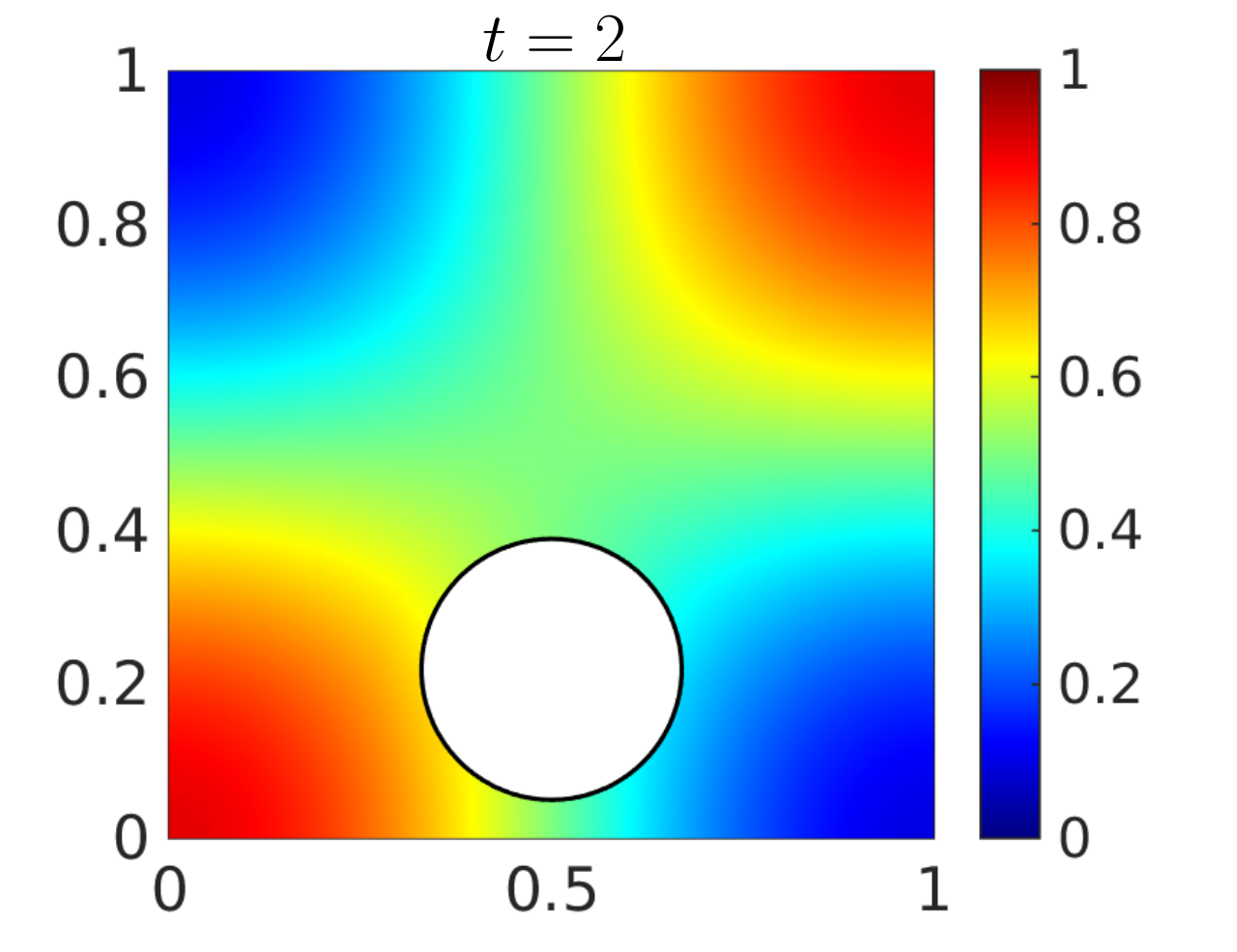}
\caption{Results for Example 1. Position of the interface and the bulk concentration on the moving interface at time t=0.5, 1.2, 2 for mesh size $h=1/40=0.025$ and time step size $k=0.5h$. \label{fig:Exexactbulkconc}}
\end{figure}

\begin{figure}\centering
\includegraphics[scale=0.4]{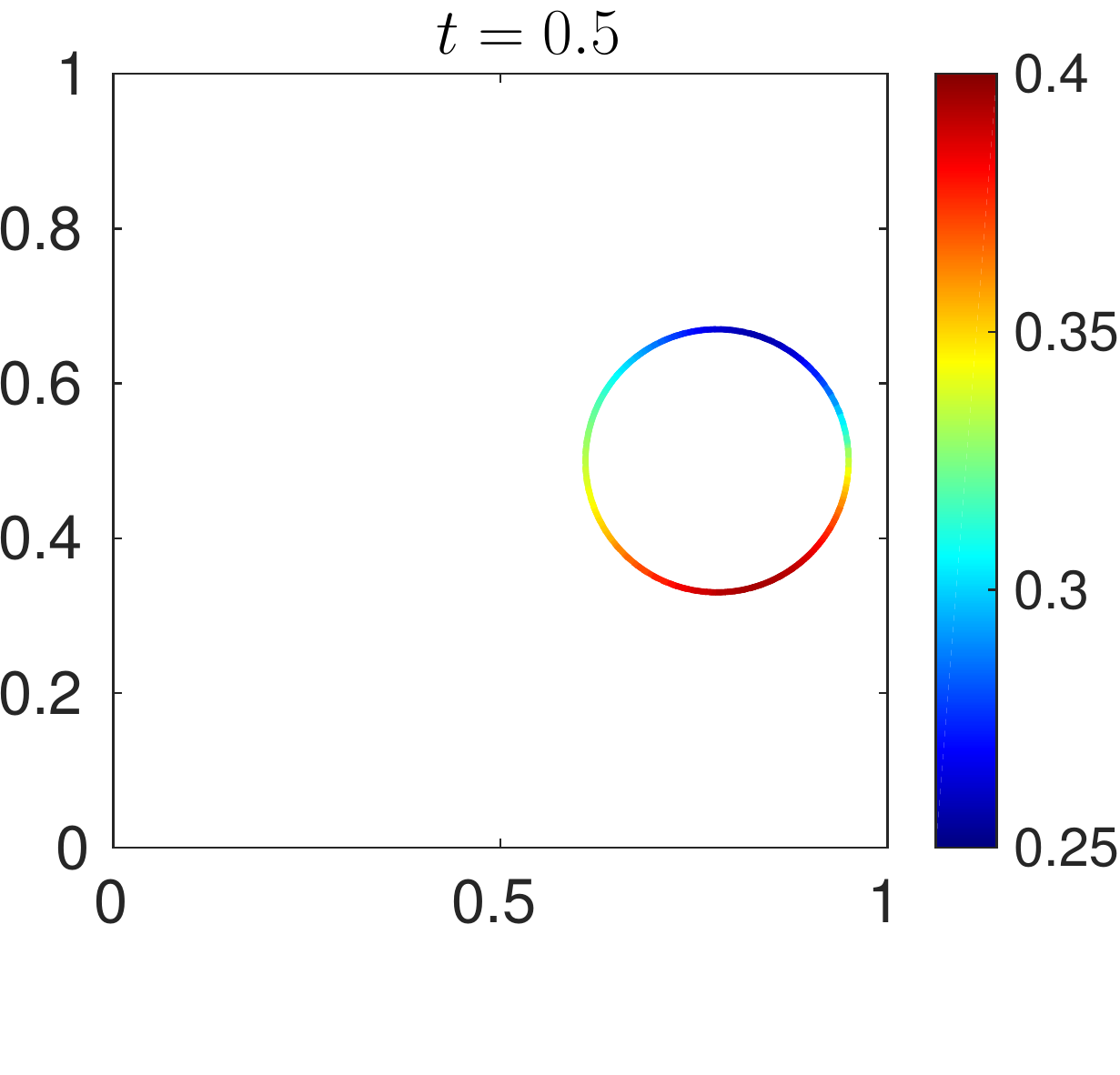} \hspace{0.1cm}
\includegraphics[scale=0.4]{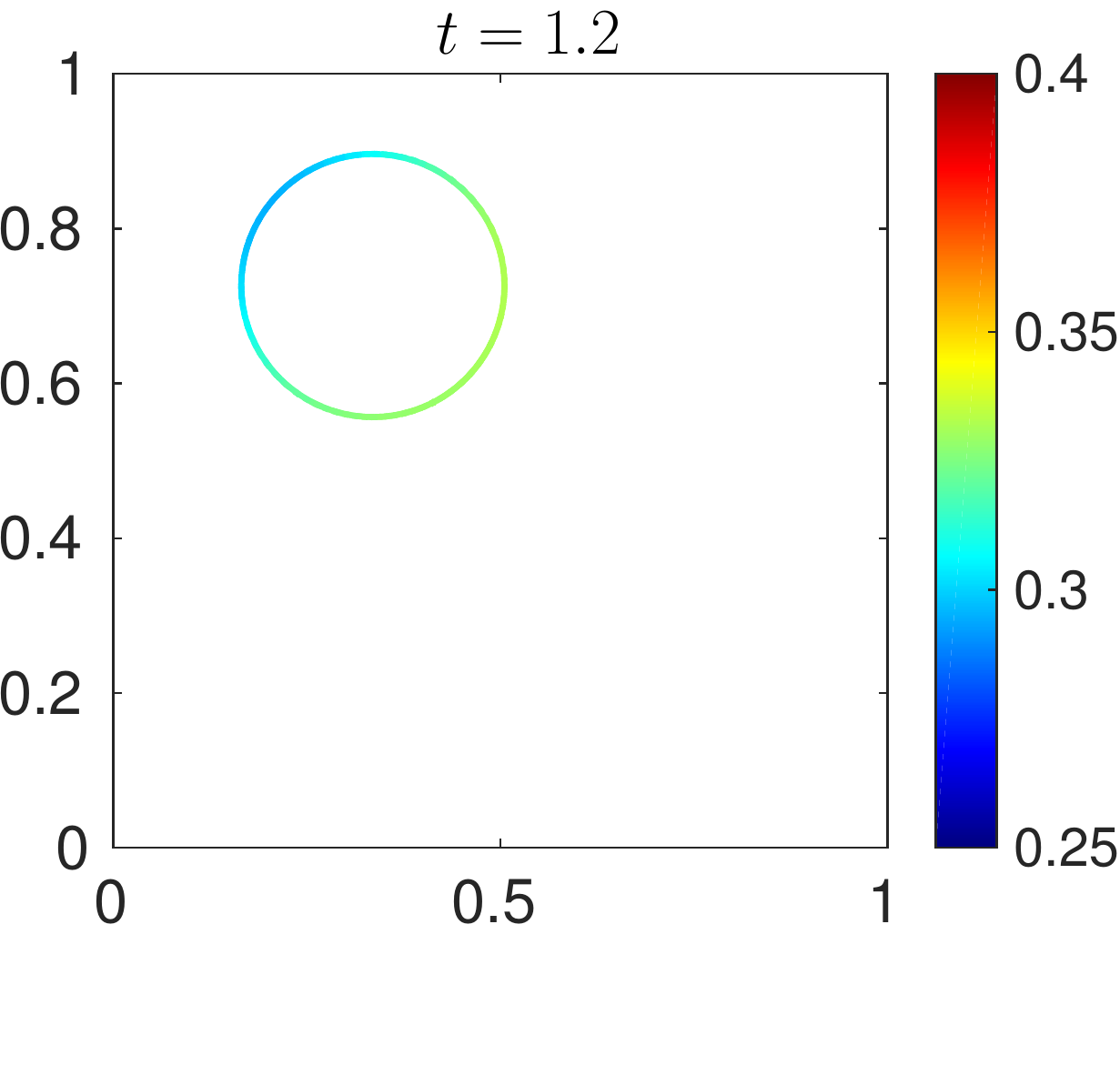}
\includegraphics[scale=0.4]{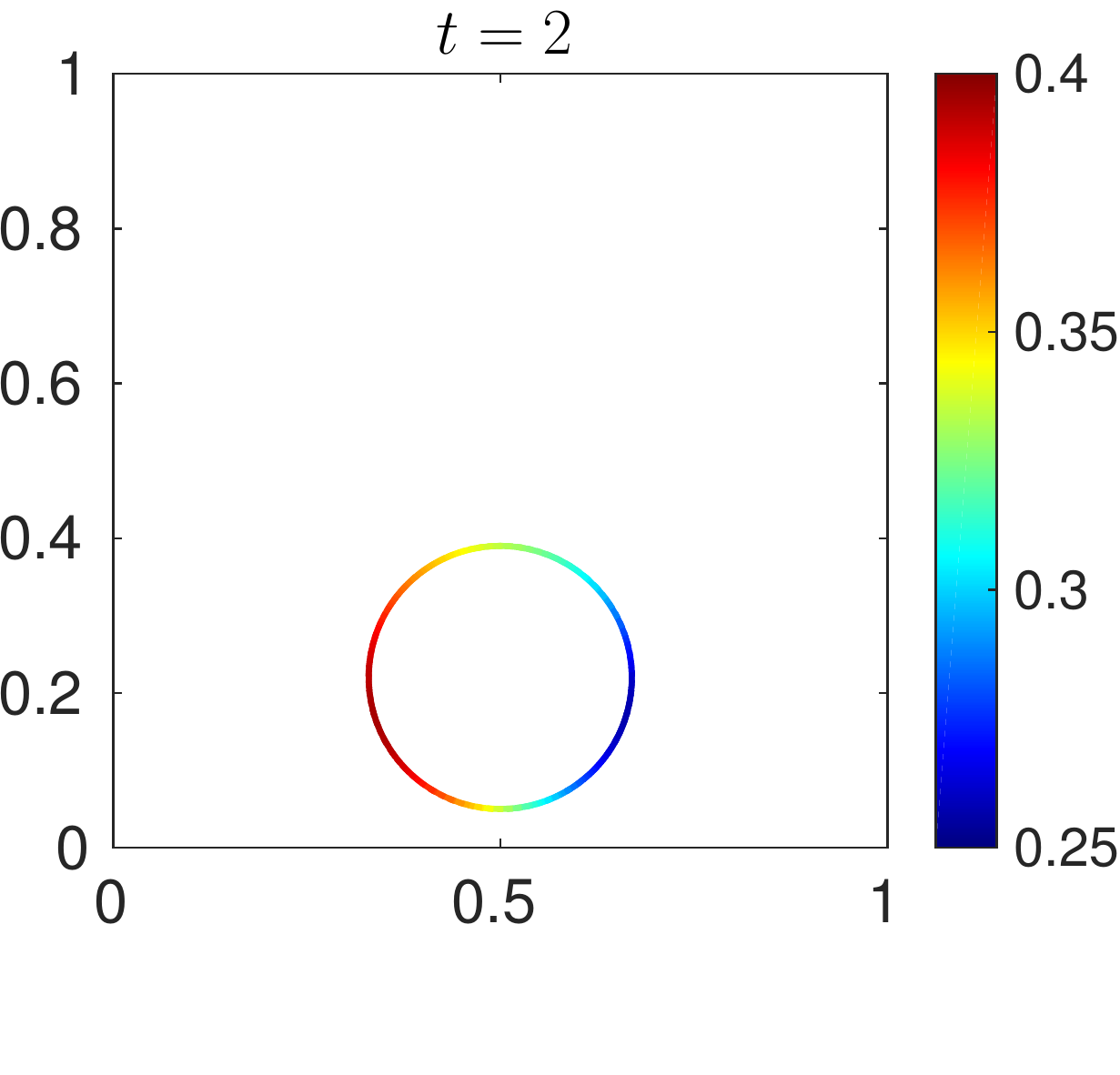}
\caption{Results for Example 1. Position of the interface and the surface concentration on the moving interface at time t=0.5, 1.2, 2 for mesh size $h=1/40=0.025$. \label{fig:Exexactsurfconc}}
\end{figure}

The bulk and interfacial concentrations on the moving interface at time $t=0.5, 1.2, 2$ are shown in Fig.~\ref{fig:Exexactbulkconc} and~\ref{fig:Exexactsurfconc}, respectively. The mesh size $h=1/40$ and the time step size $k=0.5h$.

\begin{figure}\centering
\includegraphics[scale=0.5]{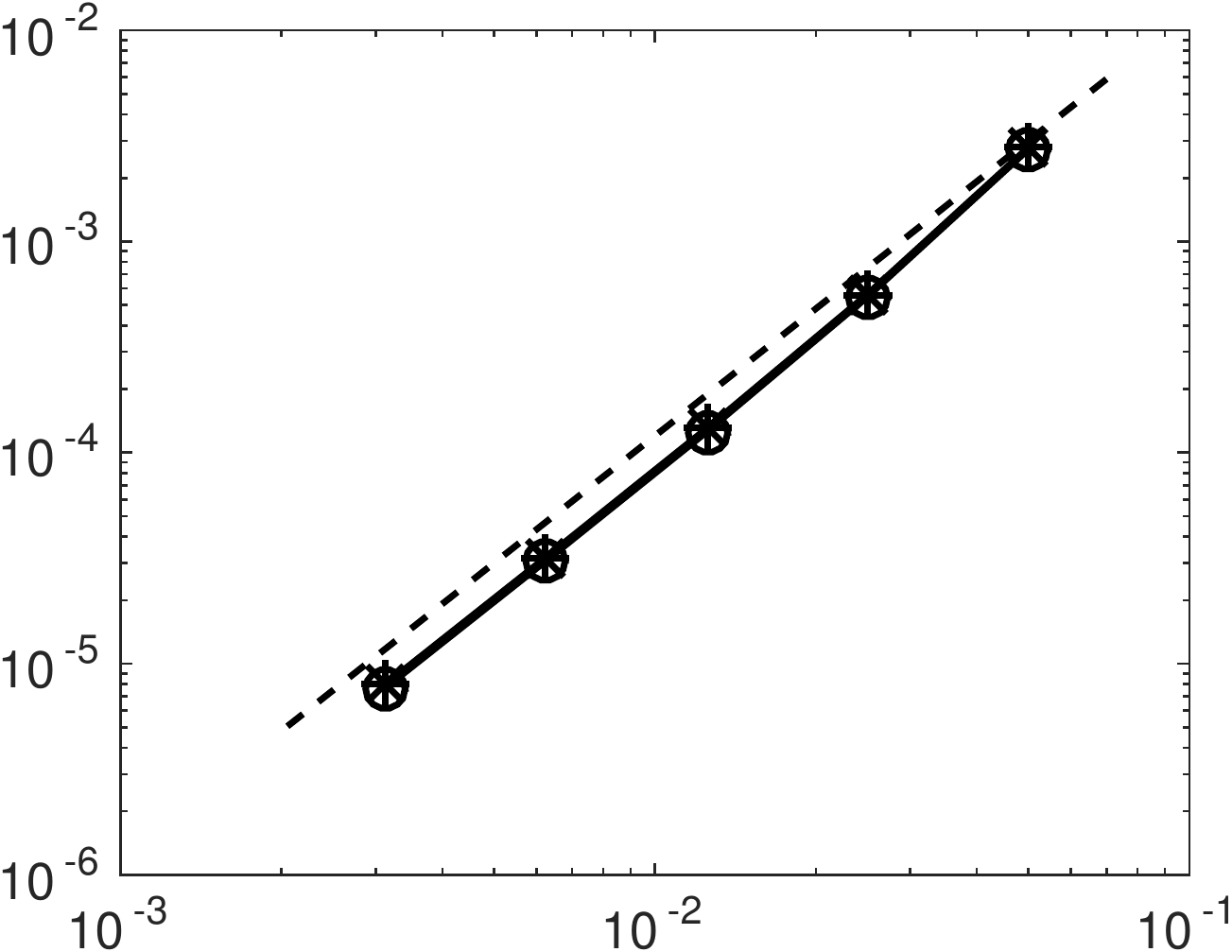}
\caption{Results for Example 1. The error $\| (u_{B} - u_{B,h}) \|_{\Omega_{h,1}(0.5)}$ measured in the $L^2$ norm versus mesh size $h$. The time step size $k=0.5h$. Results using the trapezoidal rule (stars) and the Simpson's rule (circles) for the time integration are shown. The total mass is not prescribed. The error when the total mass is prescribed coincide with the shown results. The dashed line is proportional to $h^2$. \label{fig:convuBExactsol}}
\end{figure} 

\begin{figure}\centering
\includegraphics[scale=0.5]{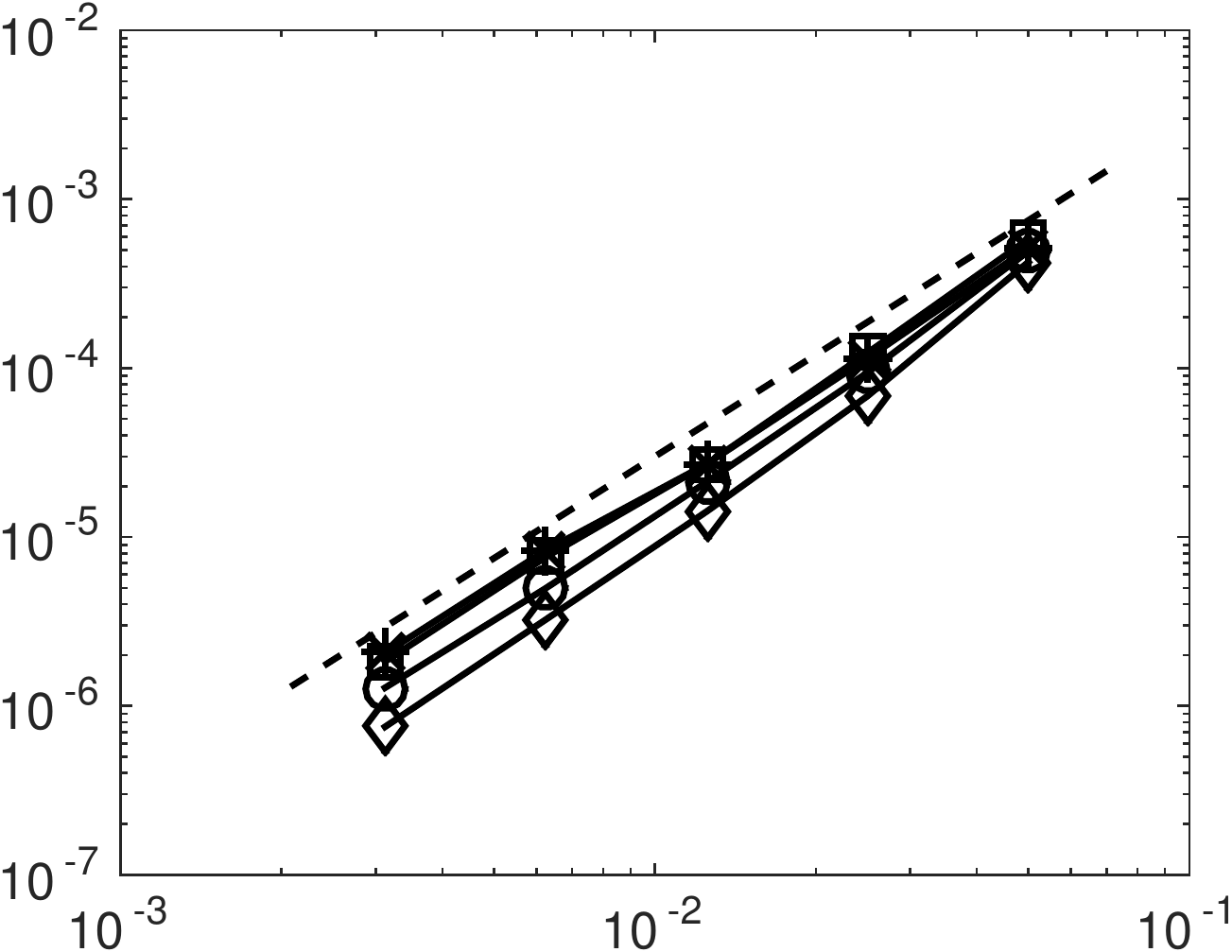}
\caption{Results for Example 1. The error $\| (u_{S} - u_{S,h}) \|_{\Gammah(0.5)}$ measured in the $L^2$ norm versus mesh size $h$. The time step size is $k=0.5h$. Squares and stars represent the trapezoidal rule with and without prescribing the total mass, respectively. Diamonds and circles represent the Simpson's rule with and without prescribing the total mass, respectively. The dashed line is proportional to $h^2$. \label{fig:convuSExactsol}}
\end{figure} 

In this example we compare the errors using the trapezoidal rule and Simpson's rule for the time integration. In Fig.~\ref{fig:convuBExactsol} and~\ref{fig:convuSExactsol} we show the errors $\| (u_{B} - u_{B,h}) \|_{\Omega_{h,1}(0.5)}$ and $\| (u_{S} - u_{S,h}) \|_{\Gammah(0.5)}$ versus mesh size $h$ using the different time integration schemes. The time step size is $k=0.5h$. We show results both with and without prescribing the total mass $\int_{\Omega_1(t)} u_B dv +\int_{\Gamma(t)} u_S ds$. In this example, the total mass is time dependent and not conserved. However, at the end of each time interval we can compute (since the exact solution is known) and prescribe the total mass. 
We observe the expected second order convergence of $\| (u_{B} - u_{B,h}) \|_{\Omega_{h,1}(0.5)}$ and $\| (u_{S} - u_{S,h}) \|_{\Gammah(0.5)}$ in the $L^2$ norm both using the trapezoidal rule and Simpson's rule with and without prescribing the total mass.  For the bulk concentration the errors with and without prescribing the total mass coincide and hence we only show the results without prescribing the total mass in Fig.~\ref{fig:convuBExactsol}. Also, the different time integration schemes give similar results because the error in the space discretization dominates. For the interfacial concentration we see in Fig~\ref{fig:convuSExactsol} that prescribing the total mass gives smaller errors than not prescribing. We also see that Simpson's rule gives smaller errors than the trapezoidal rule. Therefore, in all the other examples we use Simpson's rule and prescribe the total mass.

\begin{figure}\centering
\includegraphics[scale=0.5]{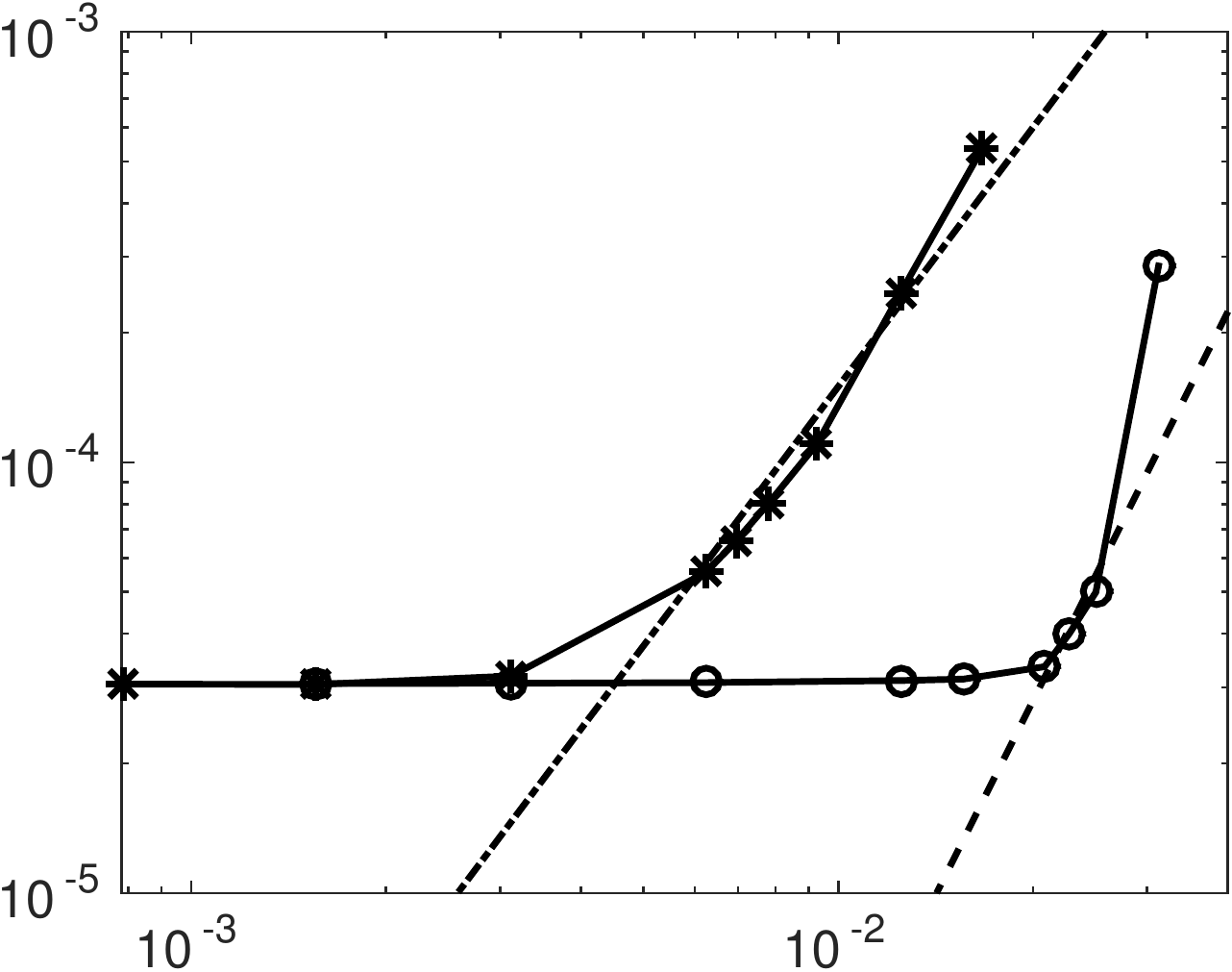}
\caption{Results for Example 1. The error $\| (u_{B} - u_{B,h}) \|_{\Omega_{h,1}(0.5)}$ measured in the $L^2$ norm versus time step size $k$ for the mesh size $h=1/160=0.00625$. Stars represent the trapezoidal rule and circles the Simpson's rule, respectively. The dashed line is proportional to $h^3$. The dashed dotted line is proportional to $h^2$ \label{fig:convuBtimeExactsol}}
\end{figure}

\begin{figure}\centering
\includegraphics[scale=0.5]{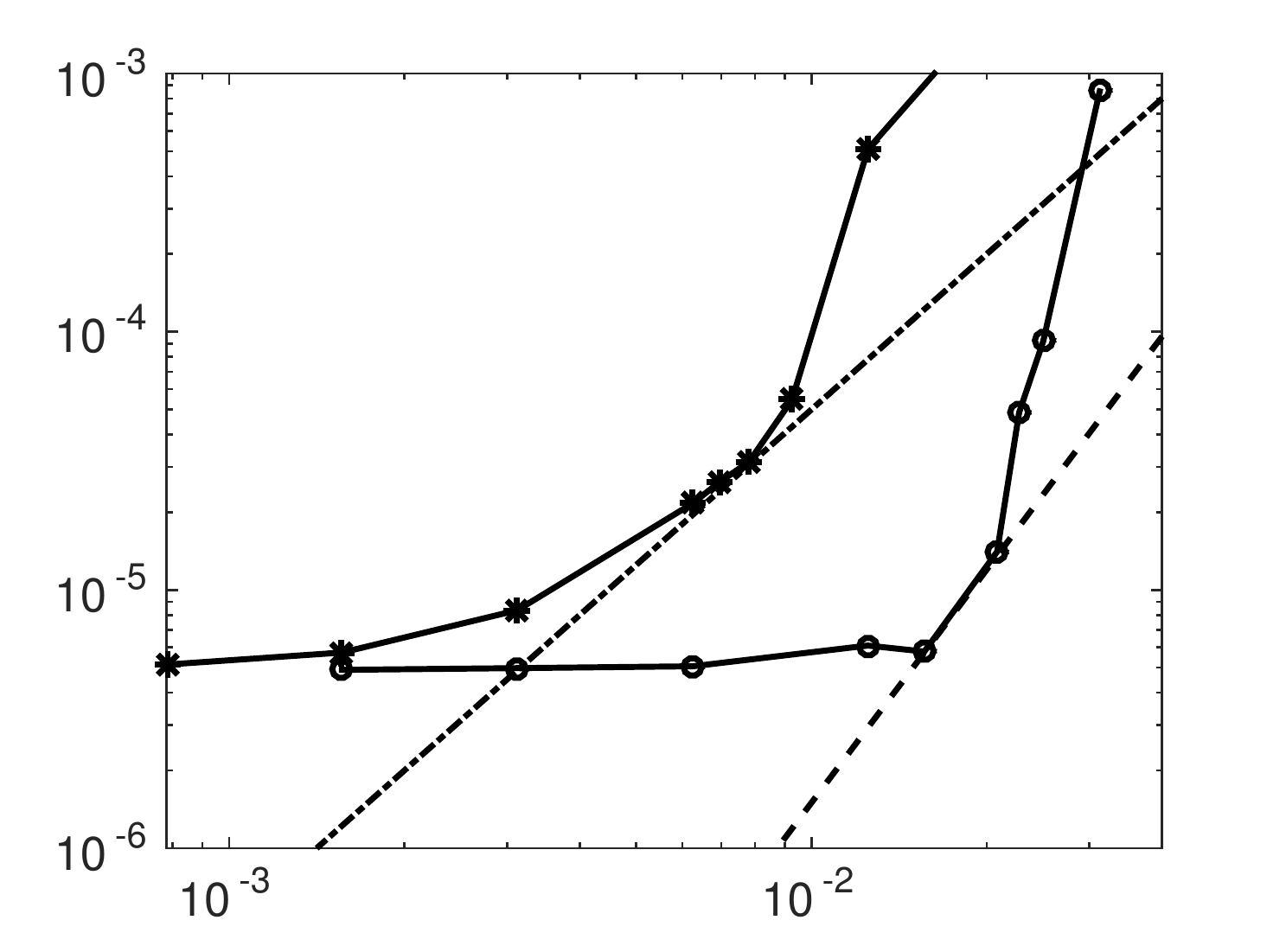}
\caption{Results for Example 1. The error $\| (u_{S} - u_{S,h}) \|_{\Gammah(0.5)}$ measured in the $L^2$ norm versus time step size $k$ for the mesh size $h=1/160=0.00625$. Stars represent the trapezoidal rule and circles the Simpson's rule, respectively. The dashed line is proportional to $h^3$. The dashed dotted line is proportional to $h^2$ \label{fig:convuStimeExactsol}}
\end{figure}

In Fig.~\ref{fig:convuBtimeExactsol} and~\ref{fig:convuStimeExactsol} we show the convergence of the bulk and interfacial concentrations at time $t=0.5$ keeping the mesh size fixed but varying the time step size. We show results both using the trapezoidal rule (stars) and Simpson's rule (circles). Using the trapezoidal rule, represented by stars in the figures, we expect to see second order convergence in time. For the bulk concentration the convergence is around second order before the space discretization dominates. For the interfacial concentration the convergence is faster for large time step sizes but as the time step decreases the error becomes around second order before the error in the space discretization dominates. For the Simpson's rule we observe the same behavior but now the order of convergence seems to be third order instead of second order. Note that the time $t=0.5$ at which the errors are measured is a nodepoint in the time discretization and since we use linear polynomials in the DG method the optimal order of convergence in time is third order, see Remark~\ref{rem:dghighorder}.

\subsection{Example 2}
\begin{figure}\centering
\includegraphics[scale=0.5]{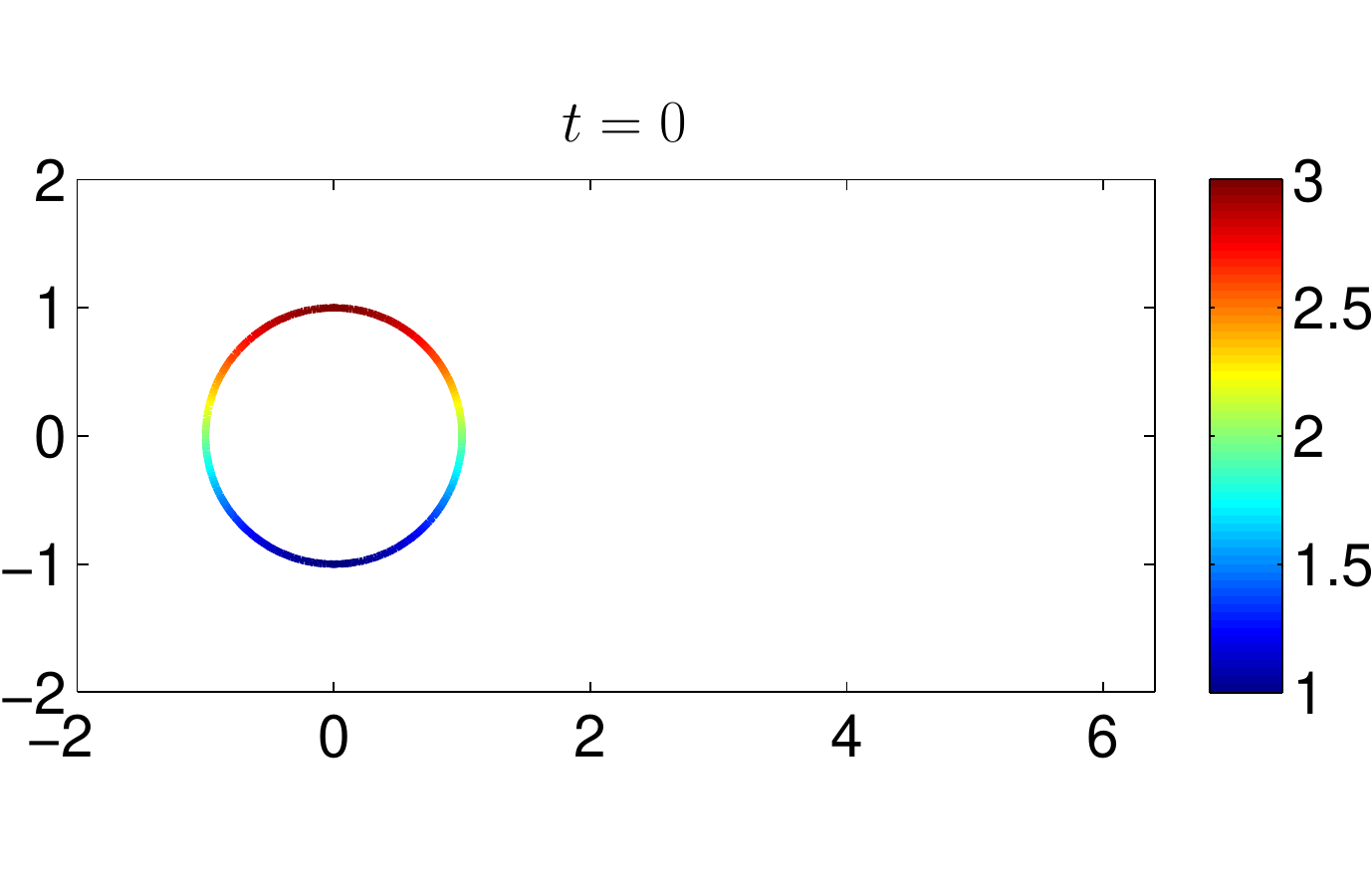} \hspace{0.1cm}
\includegraphics[scale=0.5]{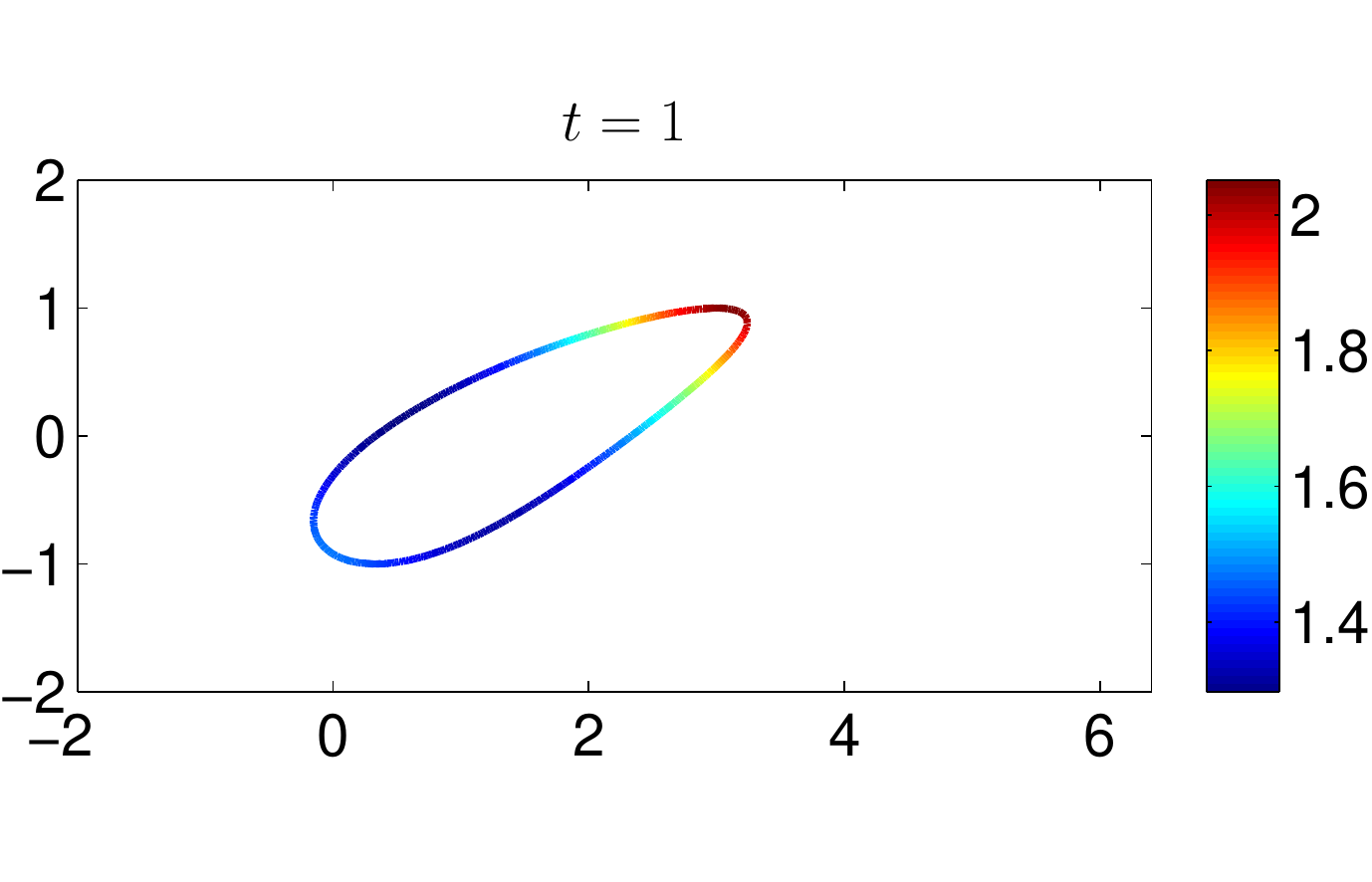}
\includegraphics[scale=0.5]{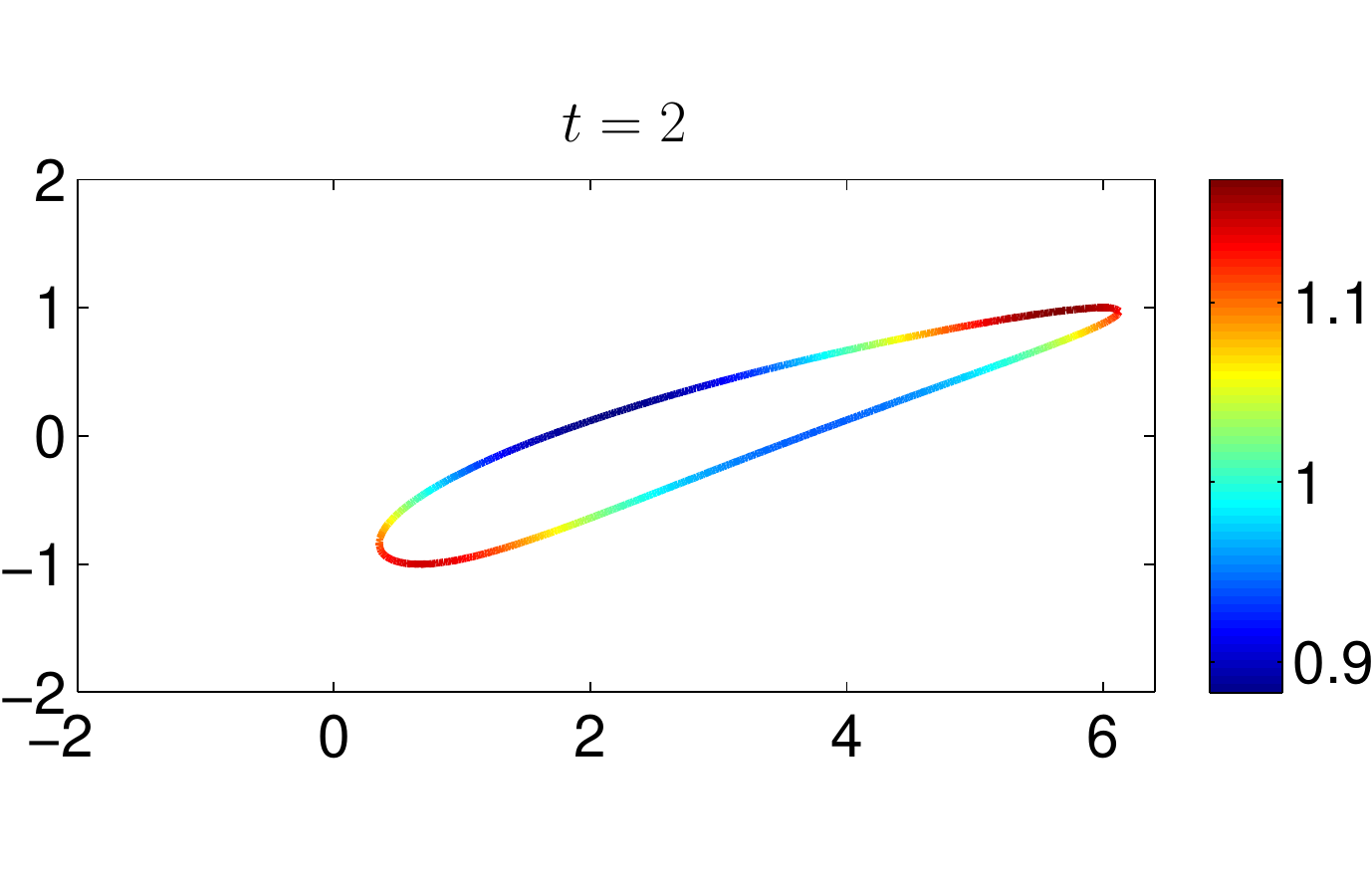}
\caption{Results for Example 2. Position of the interface and the surfactant concentration on the moving interface at time t=0, 1, 2. \label{fig:ExDIsurfconc}}
\end{figure}

\begin{figure}\centering
\includegraphics[scale=0.56]{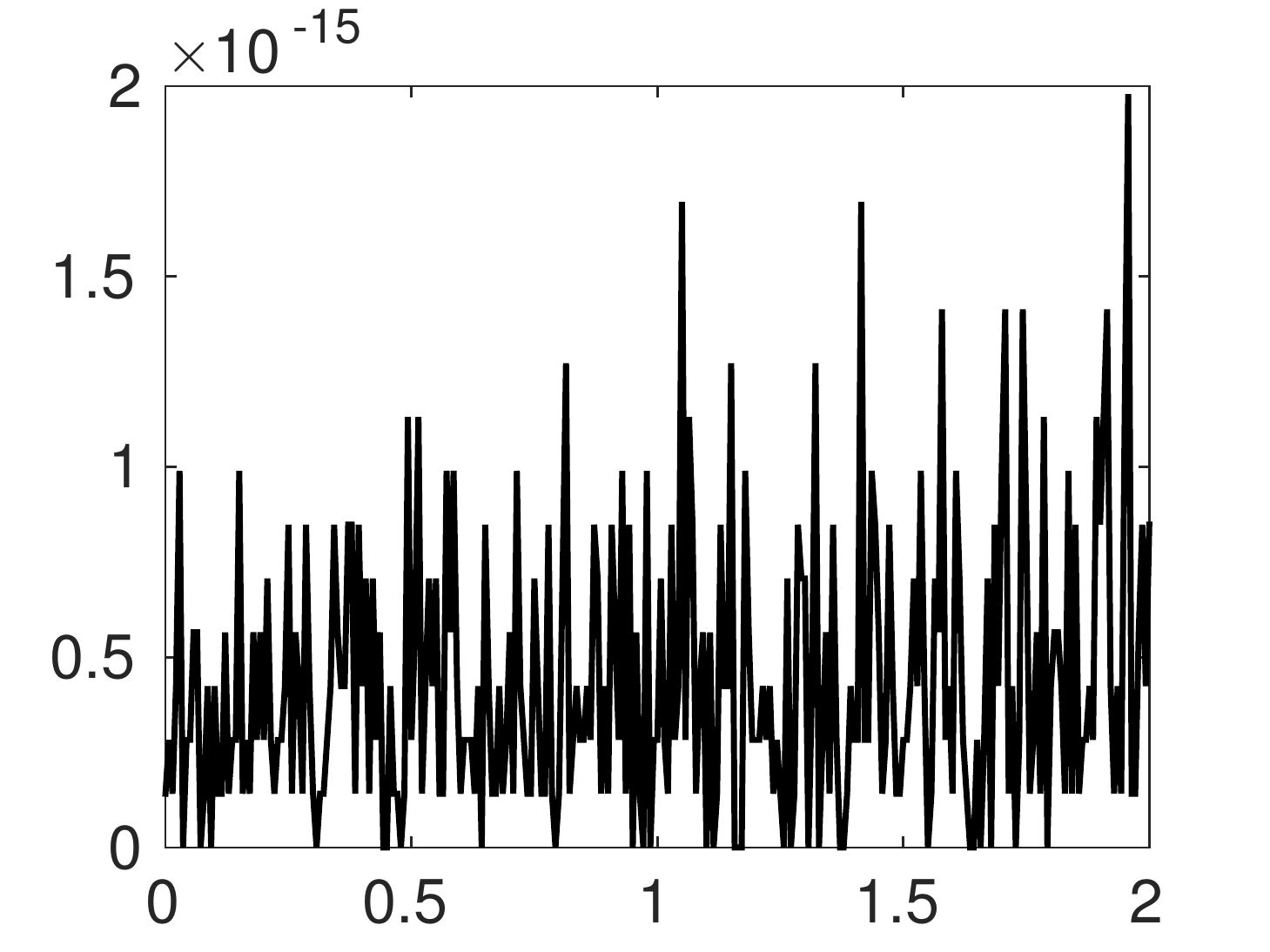}
\caption{Results for Example 2. The relative error in the total surfactant mass, $\frac{\int_{\Sigma(t)}u-4\pi r_0}{4\pi r_0}$, versus time.\label{fig:ExDIsurf}}
\end{figure}

\begin{figure}\centering
\includegraphics[scale=0.56]{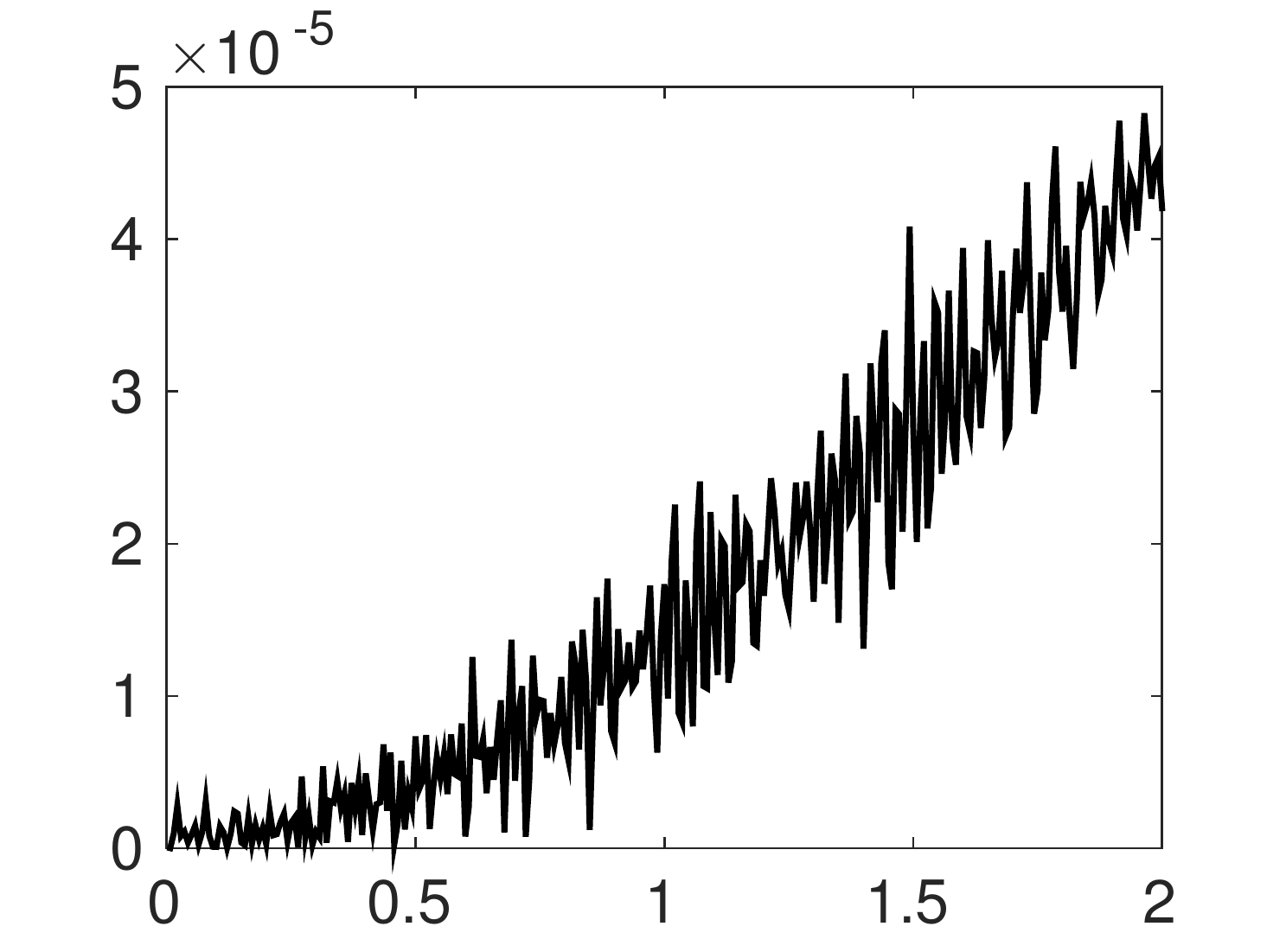}
\caption{Results for Example 2. The relative change of the area enclosed by the interface versus time.\label{fig:ExDIarea}}
\end{figure}

Next we consider Example 6 of~\cite{XuZh03}, a surface convection-diffusion problem with a moving interface modeling insoluble surfactants. Thus, surfactants only exist on the interface. Initially the interface $\intf$ is a circle centered at the origin with radius $r_0=1$ and the velocity field $\vel=(\frac{(y+2)^2}{3},0)$.  The initial interfacial surfactant concentration $\concS(0,x,y)=y/r_0+2$. The interfacial diffusion coefficient $k_S=1$. The computational domain is chosen as $\Omega=[-2, \ 6.4]\times [-2, \ 2]$ and we use 148 gridpoints in the x-direction and 71 gridpoints in the y-direction which yields a mesh size $h \approx 0.06$. The time step size $k=h/8$ as in~\cite{XuZh03}. The surfactant concentration on the moving interface at times $t=0, 1, 2$ is shown in Fig~\ref{fig:ExDIsurfconc}. In Fig.~\ref{fig:ExDIsurf} the relative error in the total surfactant mass versus time is shown and we see that the error is of the order of machine epsilon. In Fig.~\ref{fig:ExDIarea} we show the relative change of the area enclosed by the interface. We observe a change in the area by less than $0.005 \%$ at time $t=2$ for $h \approx 0.06$. In the method presented in~\cite{XuZh03} which is based on the standard level set method a surfactant mass loss of 4-5 $\%$ and a change in the area by $1\%$ for $h=0.04$ is observed, see Figs.~8 and 9 in~\cite{XuZh03}. In~\cite{XuZh03} the PDE governing the evolution of the interfacial surfactant concentration is extended off the interface to other level sets in a neighborhood of the interface (the zero level set). 

\subsection{Example 3}
We consider the same example as in Section 4.1 of~\cite{KT14}. Initially the interface $\intf$ is a circle centered in $(0,1)$ with radius $r_0=0.5$ and the velocity field $\bfbeta=(-1+y,0)$. The computational domain is chosen as $\Omega=[-1, 1] \times [0, 2]$. The non-dimensional numbers are $\textrm{Pe}_\textrm{S}=10$, $\textrm{Pe}=1$,  $\textrm{Bi}=1$, $\alpha=1$, and $\textrm{Da}=0.2$. The initial surface and bulk surfactant concentrations are $u_S(0,x,y)=0.4$ and $u_B(0,x,y)=2/3$, respectively. 

\begin{figure}
\begin{center} 
\includegraphics[width=0.45\textwidth]{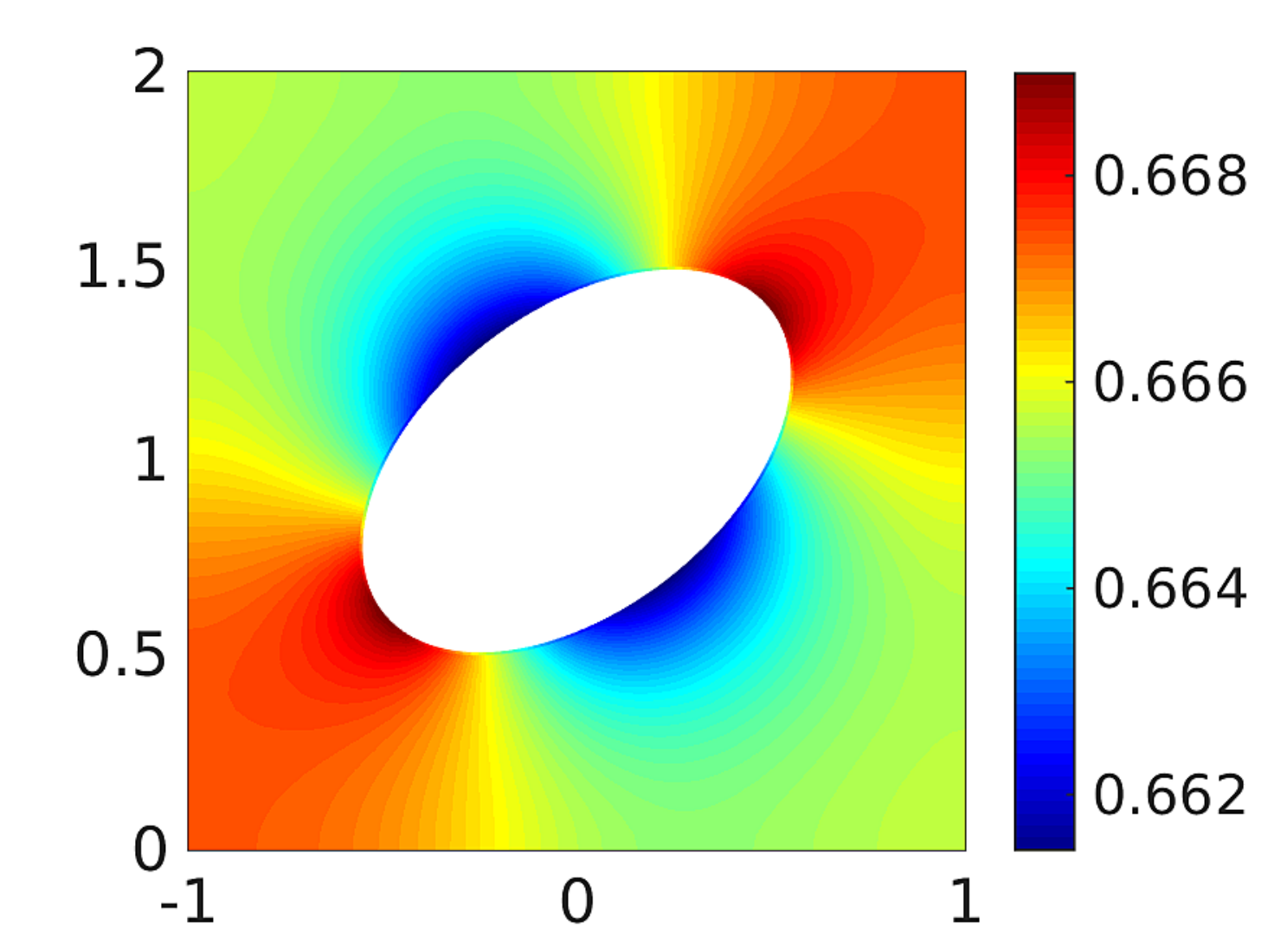}
\includegraphics[width=0.45\textwidth]{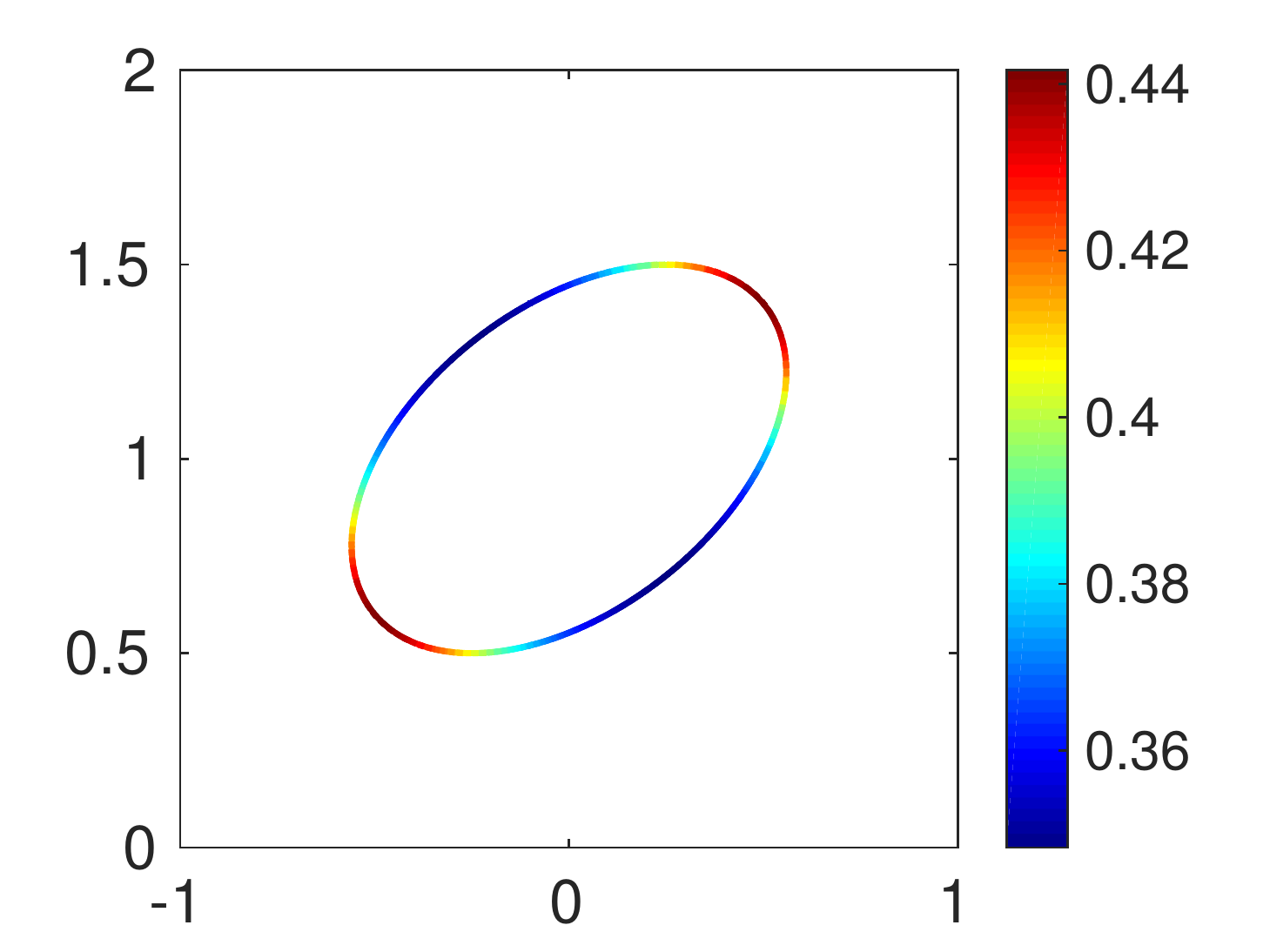}
\caption{Results for Example 3. Position of the interface and surfactant concentrations in the bulk and on the moving interface at time $t=0.5$ for mesh size $h=2/50=0.04$ and time step size $k=0.625h$.\label{fig:uBSKT}}
\end{center}
\end{figure}

\begin{figure}
\begin{center} 
\includegraphics[width=0.4\textwidth]{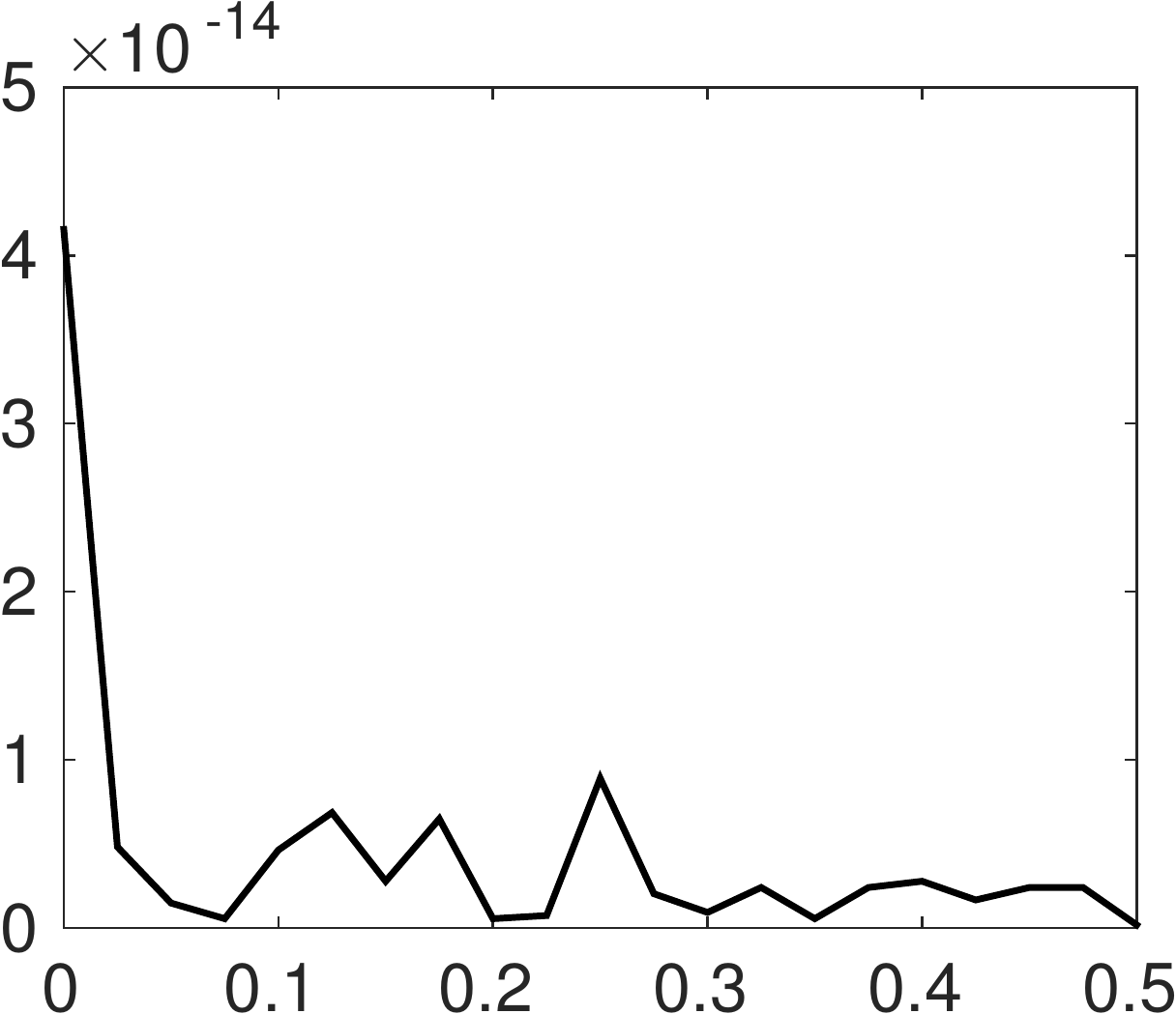}
\caption{Results for Example 3. The relative error in the total surfactant mass versus time for mesh size $h=2/50=0.04$ and time step size $k=0.625h$. \label{fig:conserv}}
\end{center}
\end{figure}

The bulk and interfacial surfactant concentrations on the moving interface at time $t=0.5$ are shown in Fig.~\ref{fig:uBSKT} for mesh size $h=2/50$ and time step size $k=0.625h$. We see in Fig.~\ref{fig:conserv} that the relative error in the total surfactant mass is of the order of machine epsilon and hence the total surfactant mass is accurately conserved. In the  method presented in~\cite{KT14} which is based on the segment projection method a surfactant mass loss of around 0.001 $\%$ is observed for the same bulk mesh as we have used but 2.25 times finer mesh for the interfacial concentration, see Fig.~8 in~\cite{KT14}.

\begin{figure}
\begin{center} 
\includegraphics[width=0.4\textwidth]{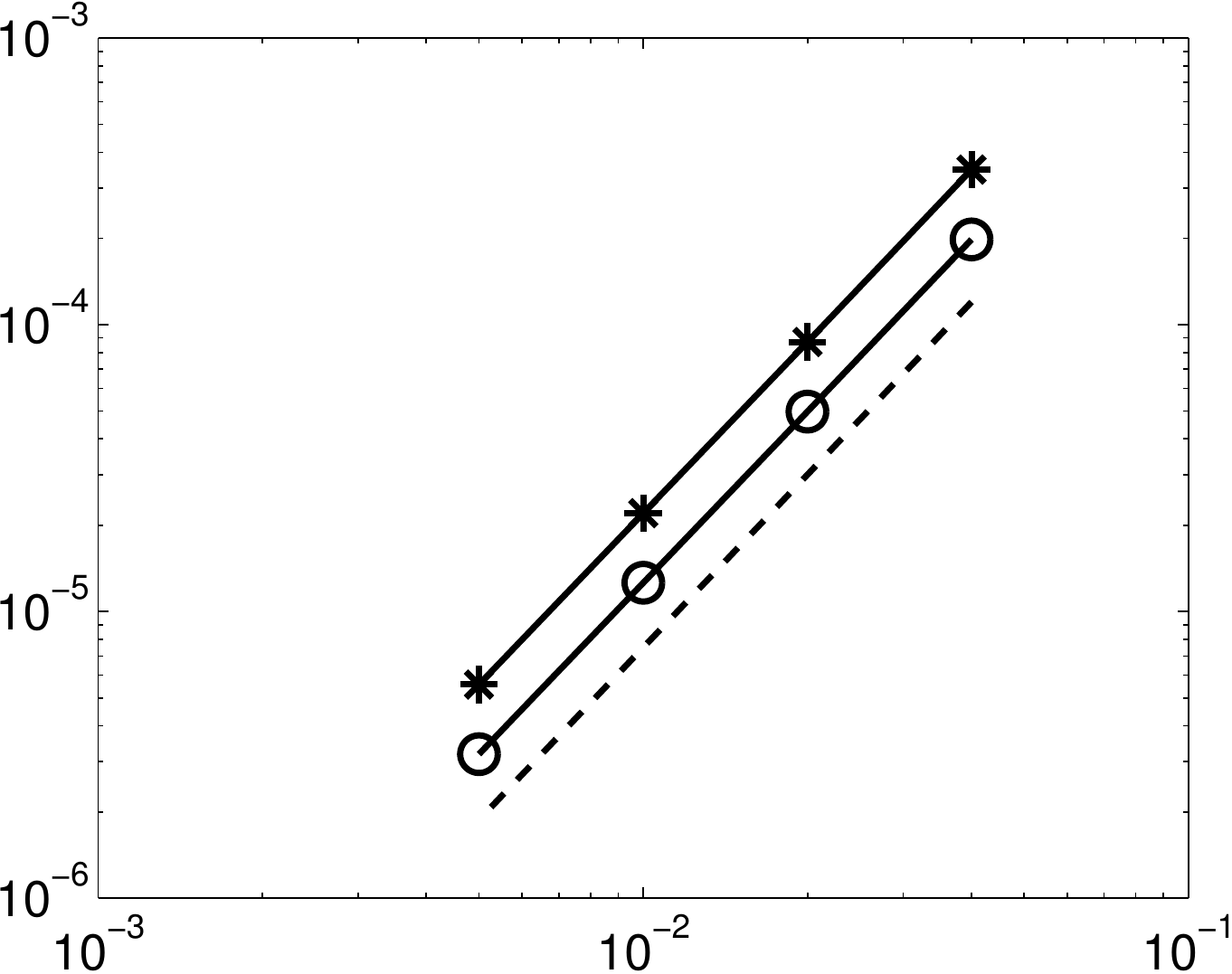}
\caption{Results for Example 3. The error $\| (u_{B} - u_{B,h}) \|_{\Omega_{h,1}(0.5)}$ measured in the $L^2$ norm (circles) and the $L^1$ norm (stars) versus mesh size $h$. The time step size $k=0.625h$.
The dashed line is proportional to $h^2$. \label{fig:convuB}}
\end{center}
\end{figure}

\begin{figure}
\begin{center} 
\includegraphics[width=0.4\textwidth]{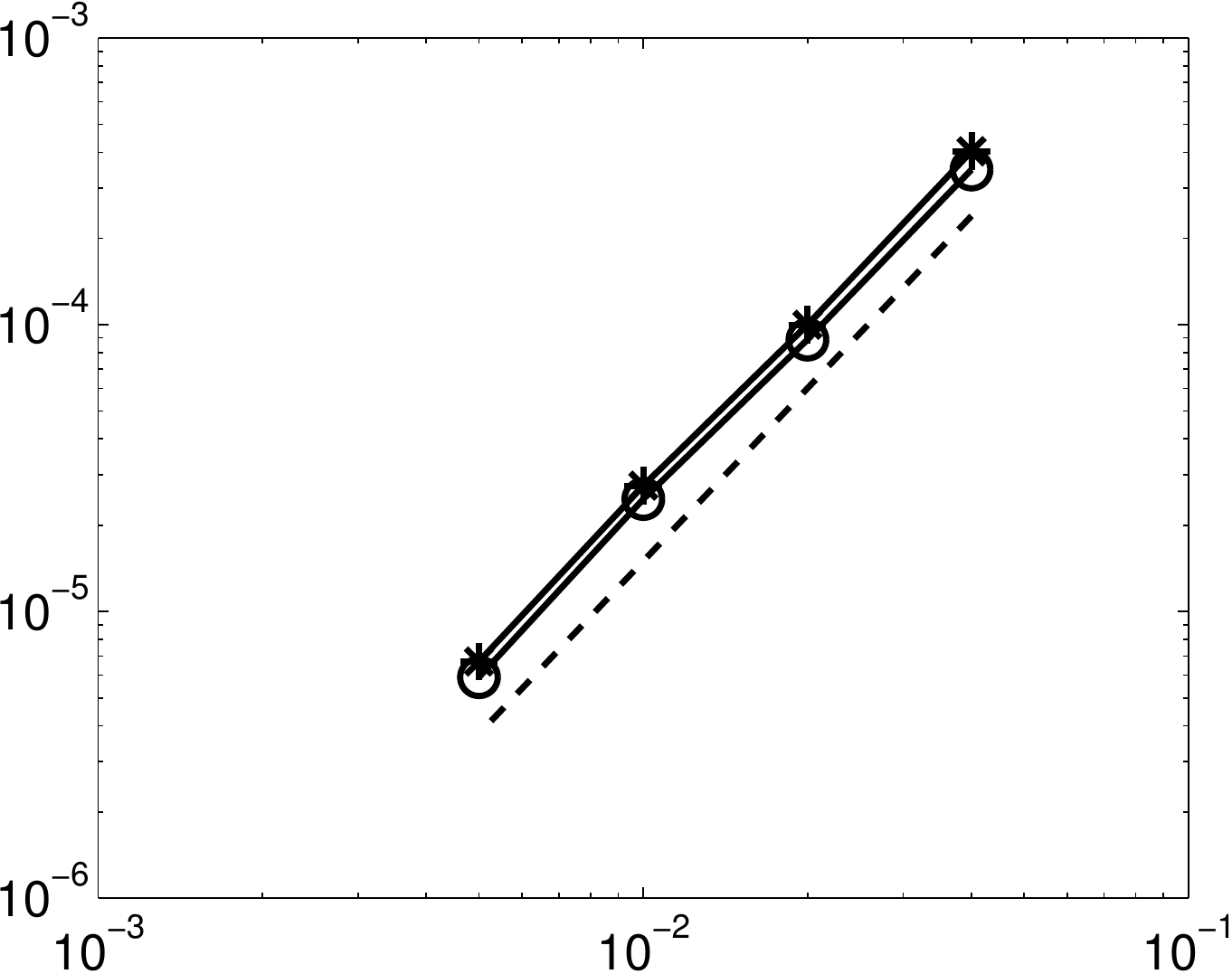}
\caption{Results for Example 3. The error $\| (u_{S,h} - u_{S,2h}) \|_{\Gammah(0.5)}$ measured in the $L^2$ norm (circles) and the $L^1$ norm (stars) versus mesh size $h$. The time step size $k=0.625h$.
The dashed line is proportional to $h^2$. \label{fig:convuS}}
\end{center}
\end{figure}

In Fig.~\ref{fig:convuB} and Fig.~\ref{fig:convuS} we show the convergence of $\| (u_{B,h} - u_{B,2h}) \|_{\Omega_{h,1}(0.5)}$ and $\| (u_{S,h} - u_{S,2h}) \|_{\Gammah(0.5)}$, respectively. We observe second order convergence both in the $L^2$ norm and the $L^1$ norm. We use the same mesh sizes as used in~\cite{KT14} for the bulk concentration. The method in~\cite{KT14} is also second order. It uses Strang splitting to both split the coupled bulk-surface problem and to split the advection and the diffusion parts of the convection diffusion equations. In the discretization of the bulk PDE regular finite difference stencils are used which require points outside of the domain $\Omega_1$. These values are determined such that the boundary conditions are enforced using an embedded boundary method. Errors measured in the $L^1$ norm reported in Fig. 8 of~\cite{KT14} are smaller than we observe. This could be explained by the fact that in~\cite{KT14} the interfacial surfactant is discretized on 2.25 times finer meshes than the meshes we use and the approximation of the interface in~\cite{KT14} is more accurate than the approximation we have used. However, in comparison the presented method is very simple to implement both in two and three space dimensions.

\subsection{Example 4}
We consider here the same coupled bulk-surface problem as in Section 5.3 of~\cite{ChLai14}. The initial interface is a circle with radius $r_0=0.3$ centered in $(0.1,0)$ and the velocity field is given by 
\begin{equation}
\bfbeta=\left(-\frac{1}{2}(1+\cos(\pi x))\sin(\pi y), \frac{1}{2}(1+\cos(\pi y))\sin(\pi x) \right)
\end{equation} 
The computational domain is chosen as $\Omega=[-1, 1] \times [-1, 1]$. The non-dimensional numbers are set to $\textrm{Pe}=\textrm{Pe}_\textrm{S}=100$ and $\textrm{Bi}=\alpha=\textrm{Da}=1$. The initial surface and bulk surfactant concentrations are $u_S(0,x,y)=0$ and 
\begin{equation}
u_B(0,x,y)=\left\{
\begin{array}{l}
0.5(1-x^2)^2 \quad \textrm{if $r>1.5r_0$} \nonumber \\
0.5(1-x^2)^2w(r) \quad \textrm{if $r_0\leq r \leq 1.5r_0$} \nonumber \\
0 \quad \textrm{otherwise} 
\end{array}\right.
\end{equation}
with $r=\sqrt{(x-x_0)^2+(y-y_0)^2}$ and
\begin{equation}
w(r)=\frac{1}{2}\left(1-\cos \left(\frac{(r-r_0)\pi}{0.5r_0}\right)\right)
\end{equation}

\begin{figure}
\begin{center} 
\includegraphics[width=0.45\textwidth]{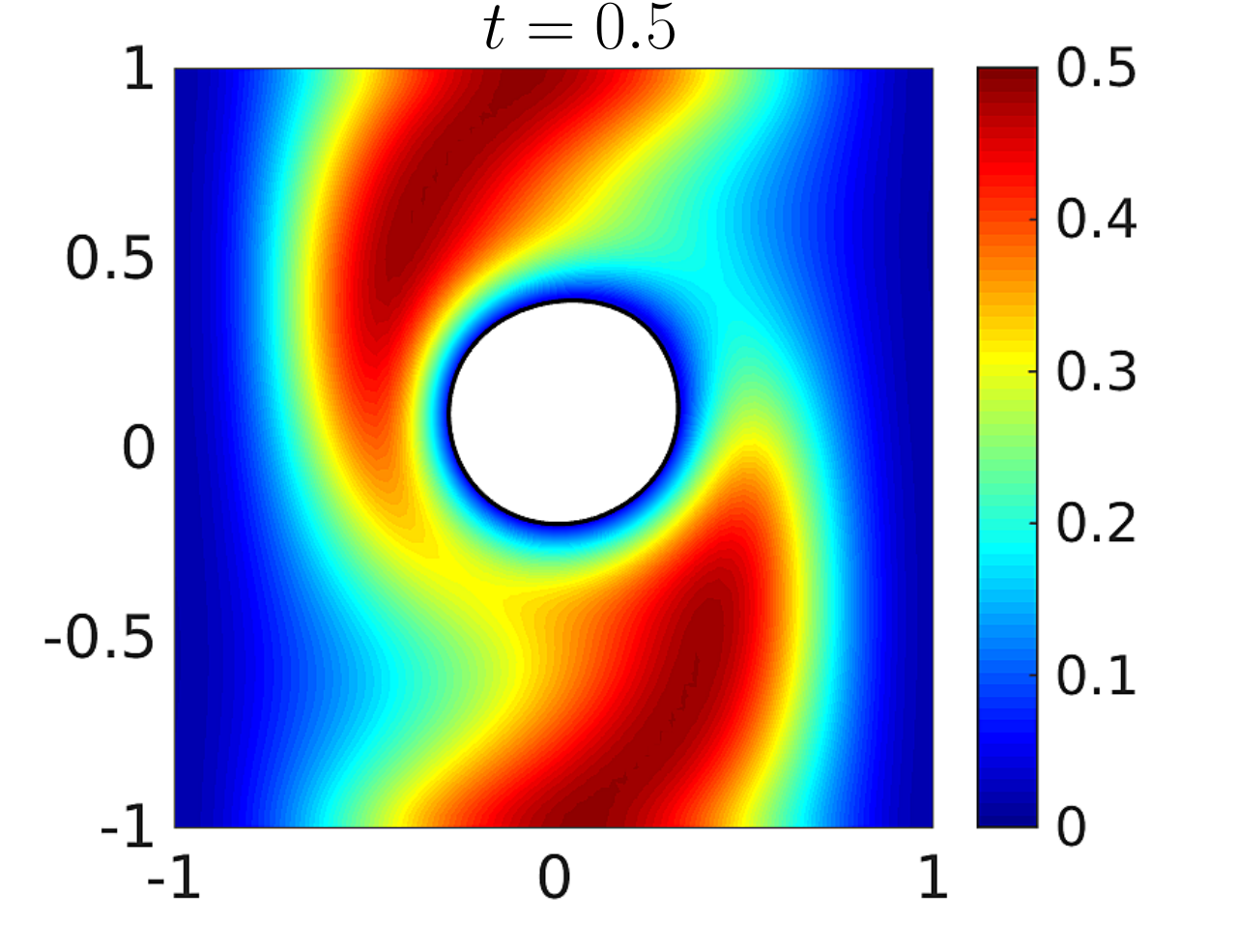}
\includegraphics[width=0.45\textwidth]{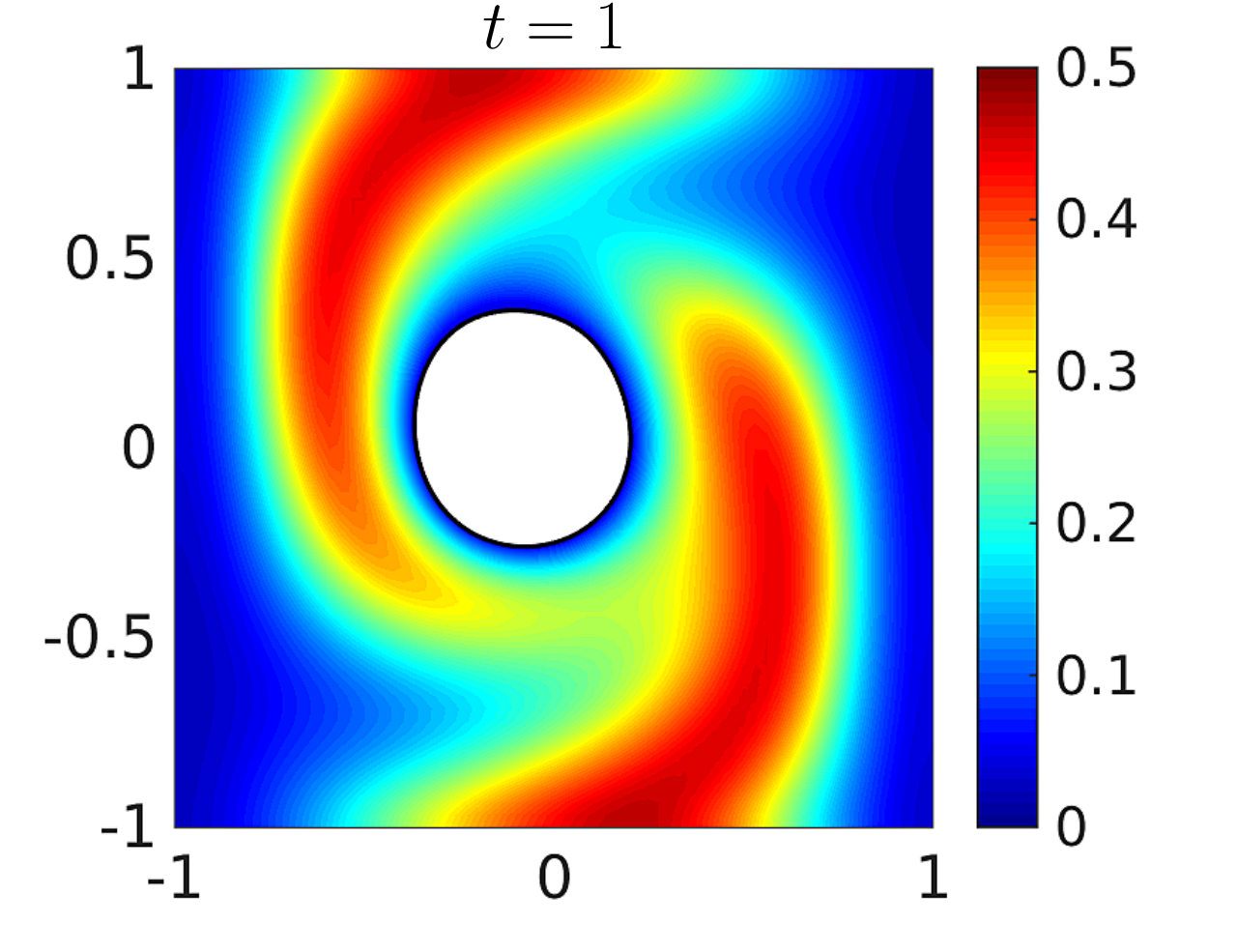}\\
\includegraphics[width=0.45\textwidth]{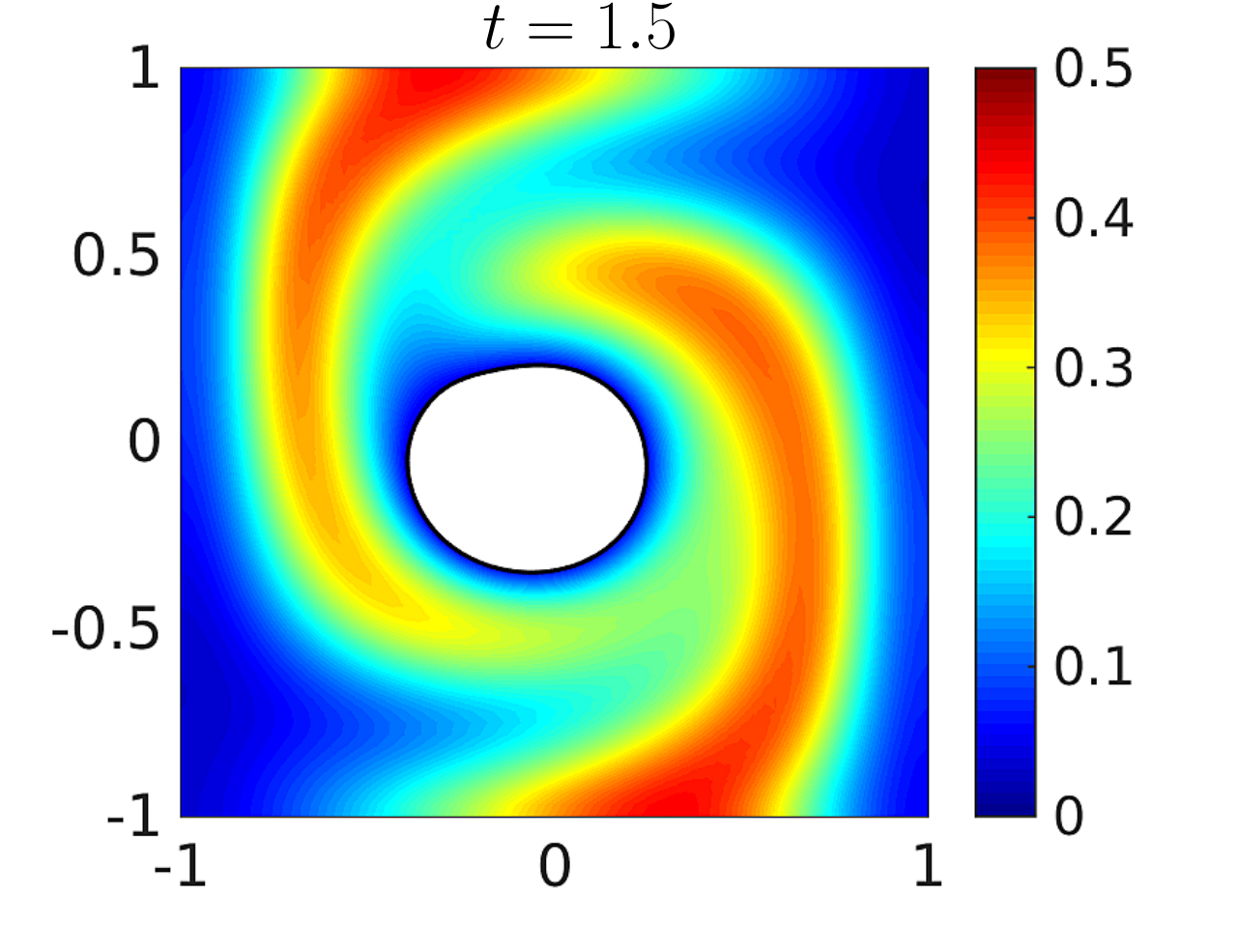}
\includegraphics[width=0.45\textwidth]{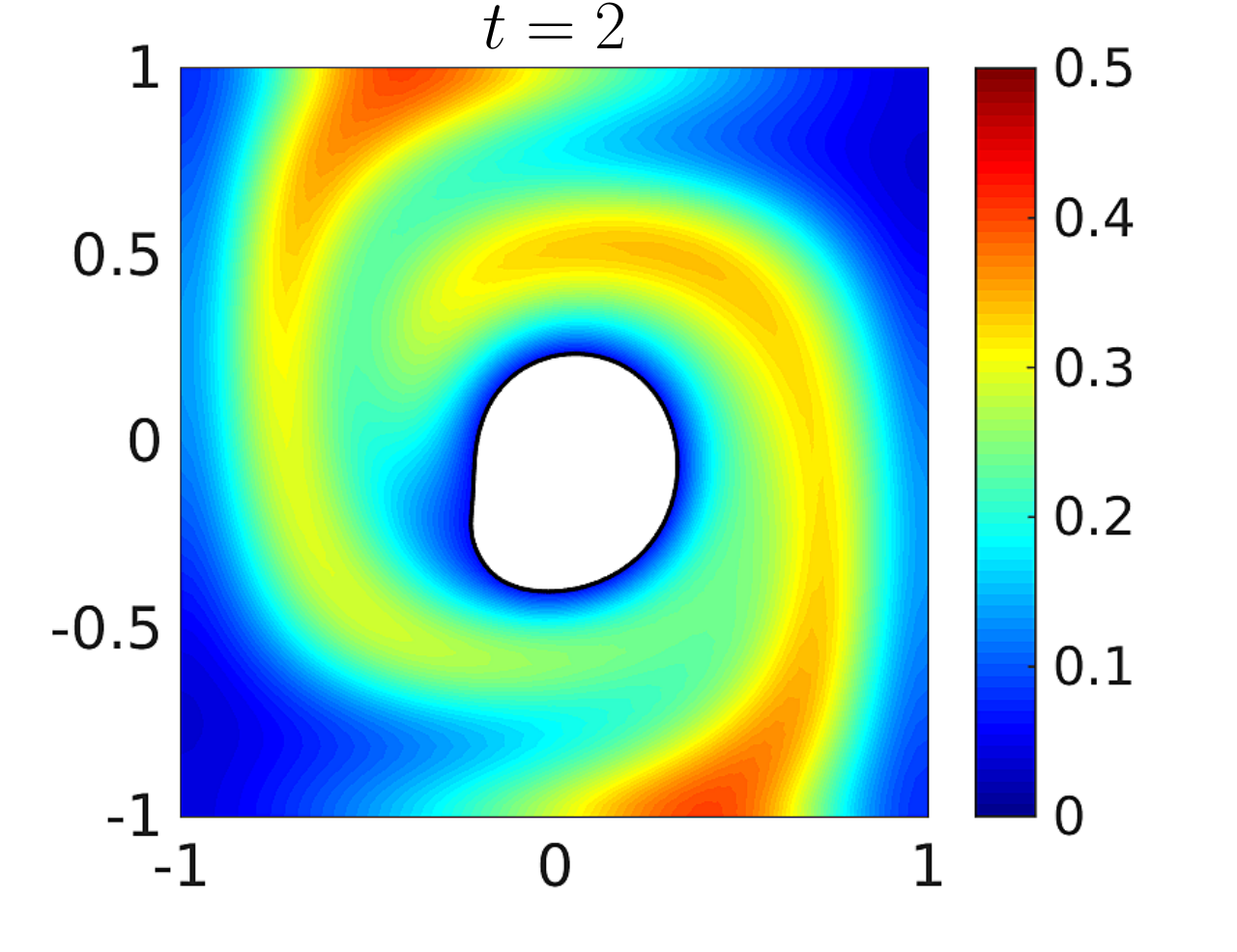}
\caption{Results for Example 4. Position of the interface and the bulk concentration at time $t=0.5, 1, 1.5, 2$ for mesh size $h=2/64=0.03125$ and time step size $k=h/8$.\label{fig:uBLai}}
\end{center}
\end{figure}

\begin{figure}
\begin{center} 
\includegraphics[width=0.45\textwidth]{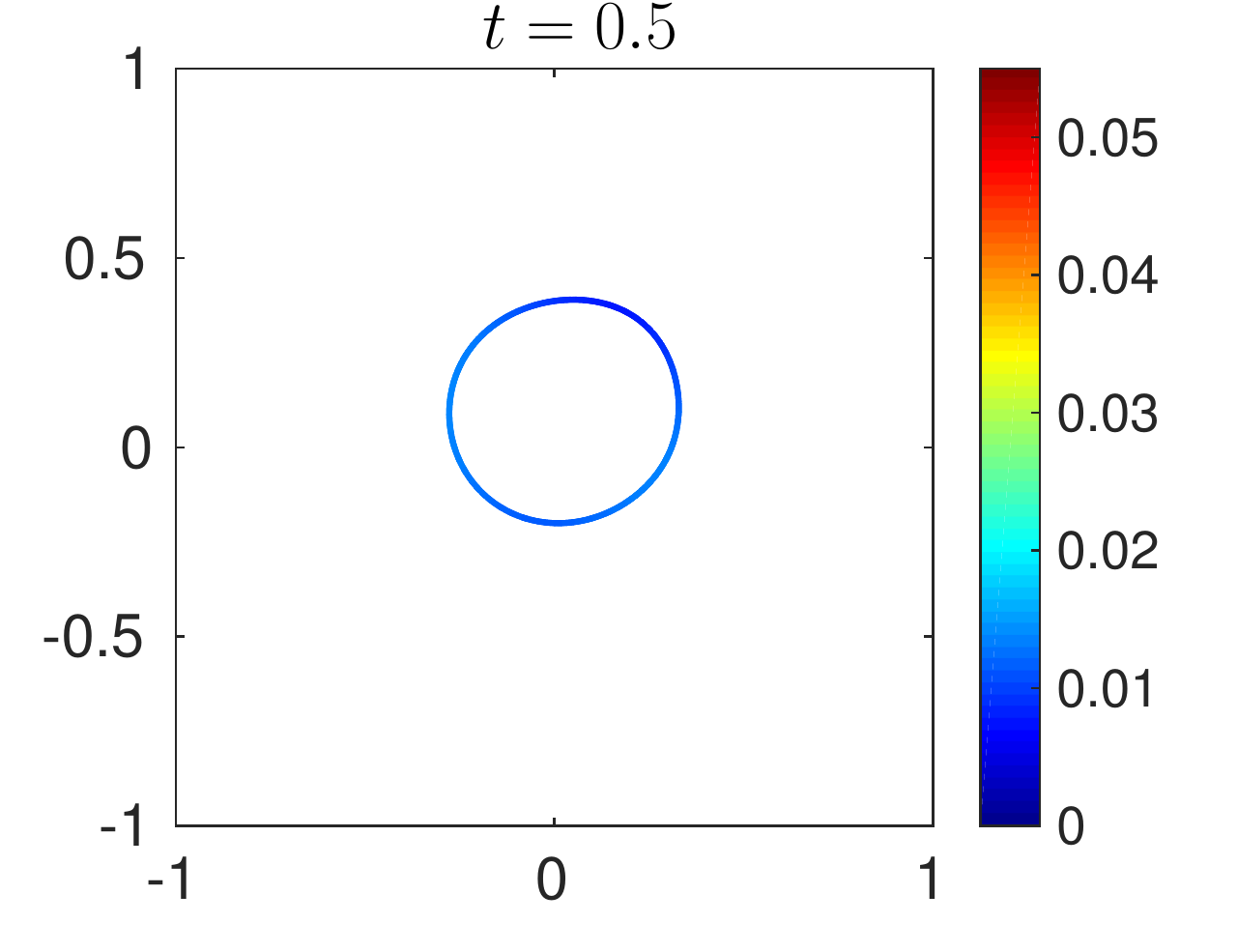}
\includegraphics[width=0.45\textwidth]{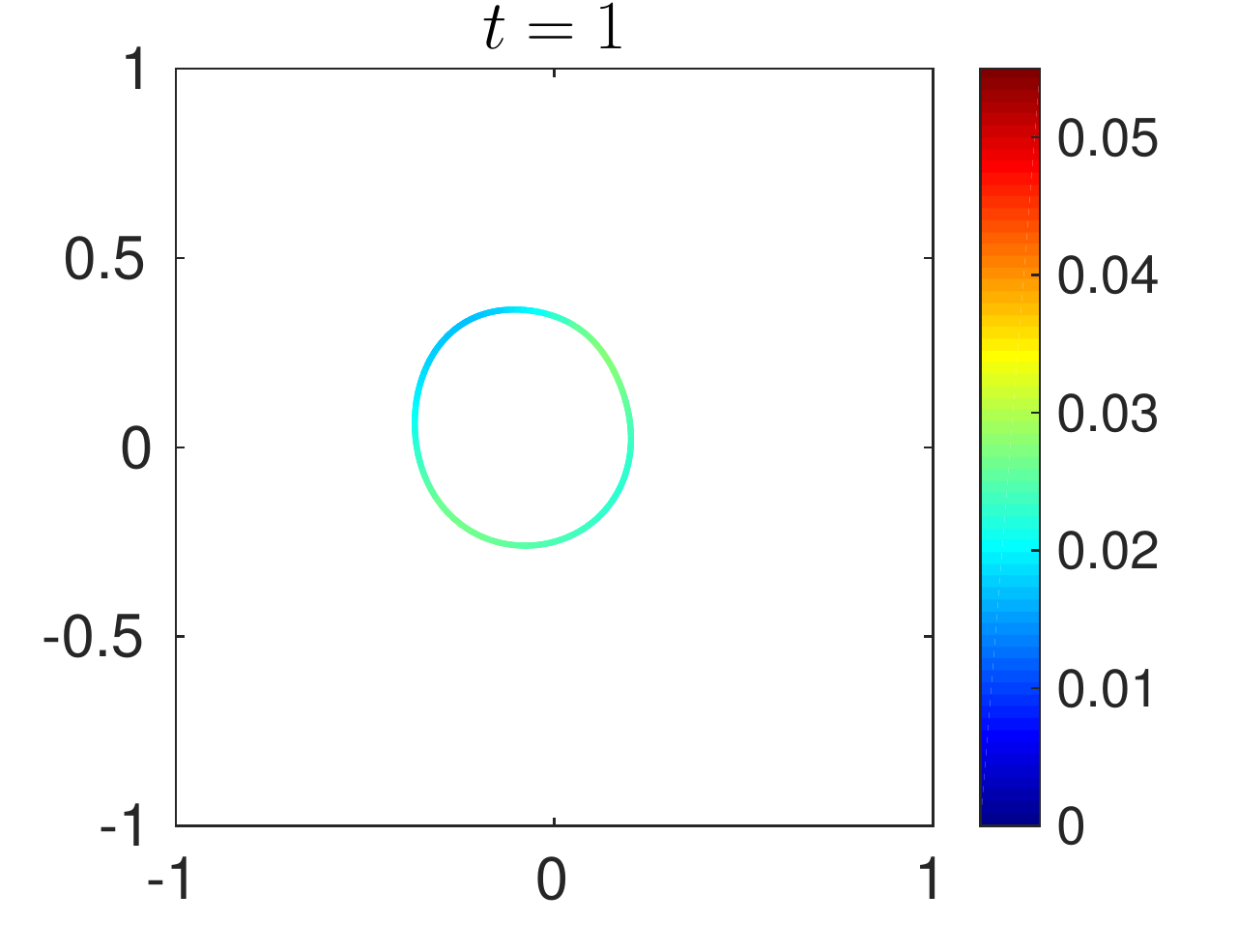}\\
\includegraphics[width=0.45\textwidth]{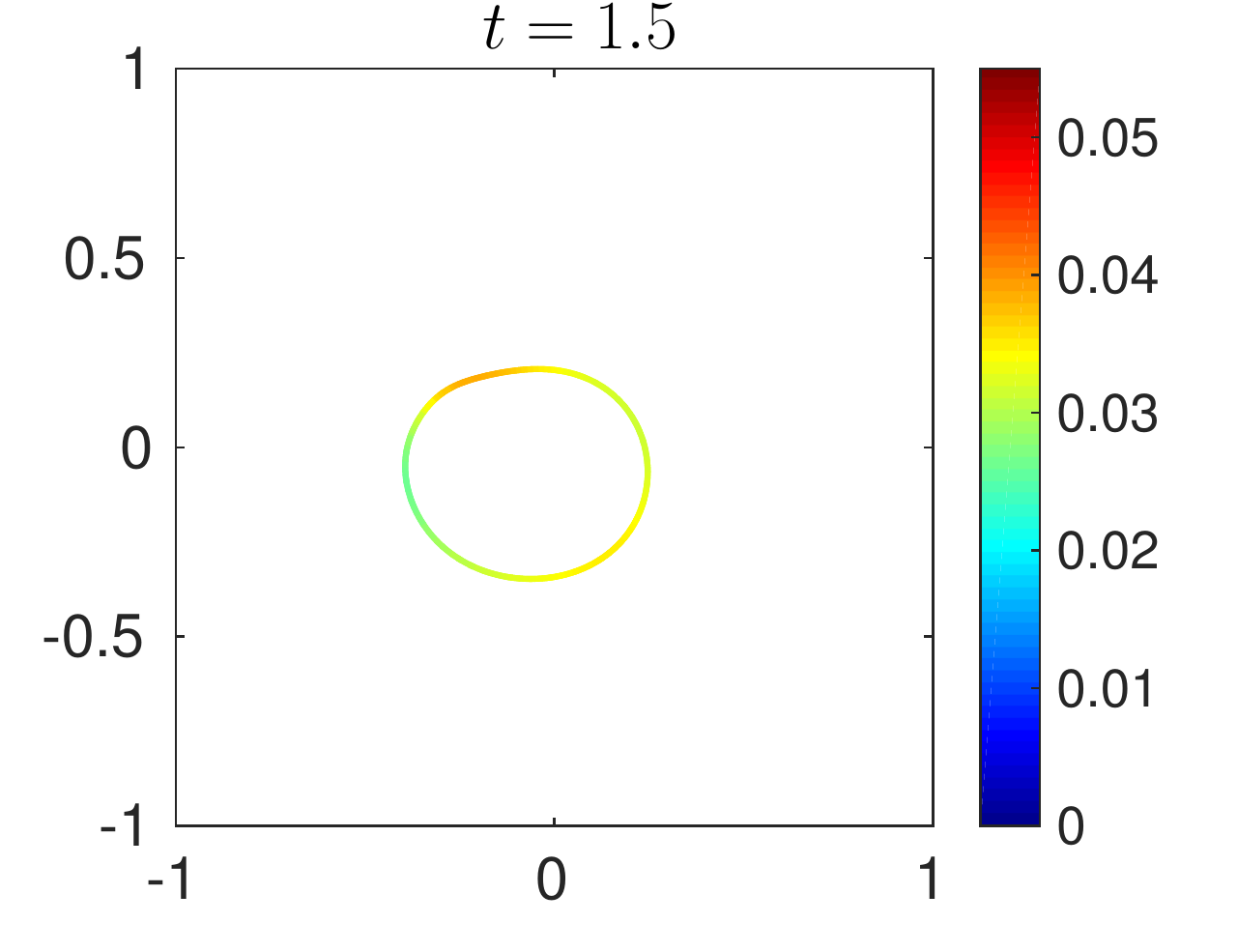}
\includegraphics[width=0.45\textwidth]{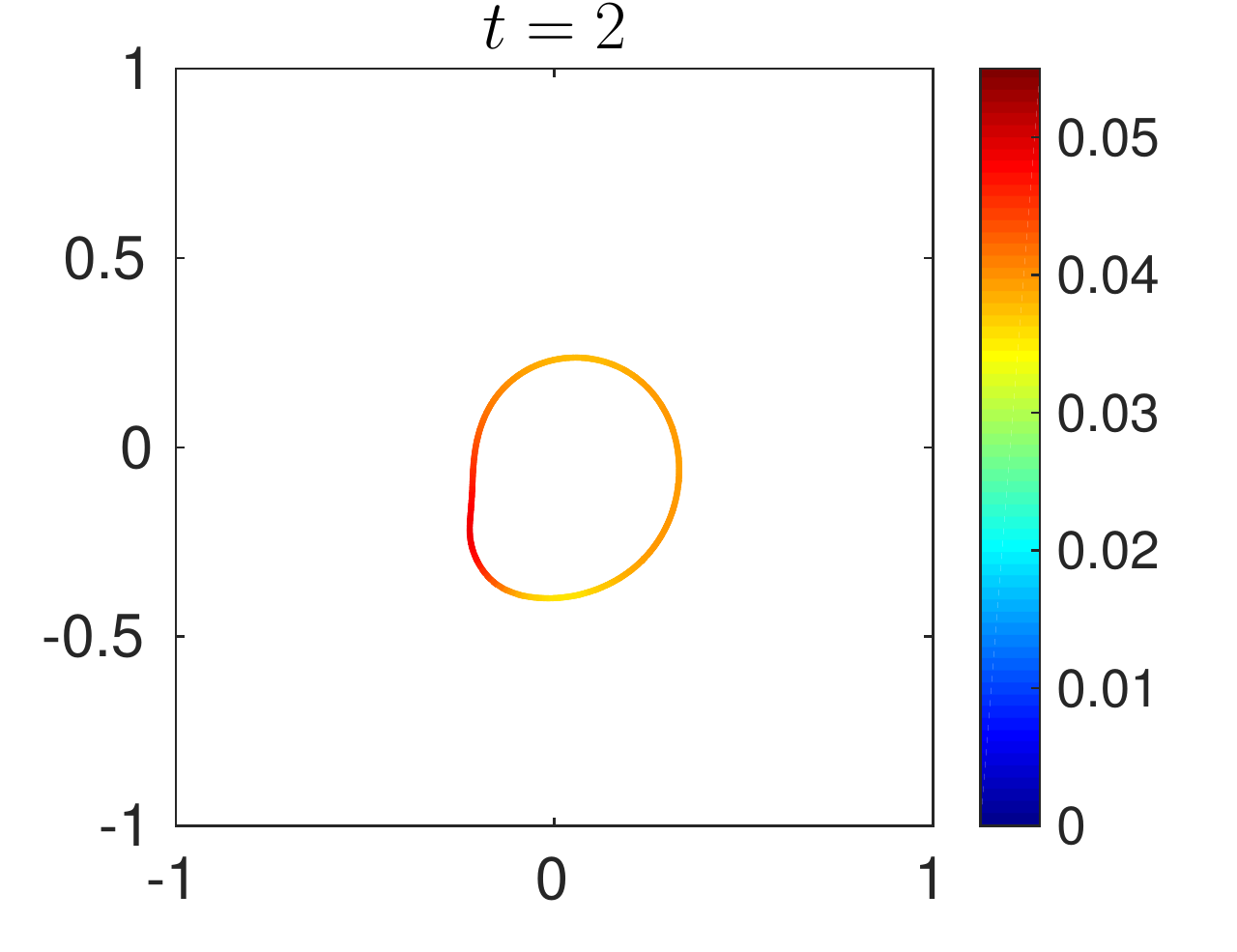}
\caption{Results for Example 4. Position of the interface and the surface concentration at time $t=0.5, 1, 1.5, 2$ for mesh size $h=2/64=0.03125$ and time step size $k=h/8$. \label{fig:uSLai}}
\end{center}
\end{figure}

\begin{figure}
\begin{center} 
\includegraphics[width=0.4\textwidth]{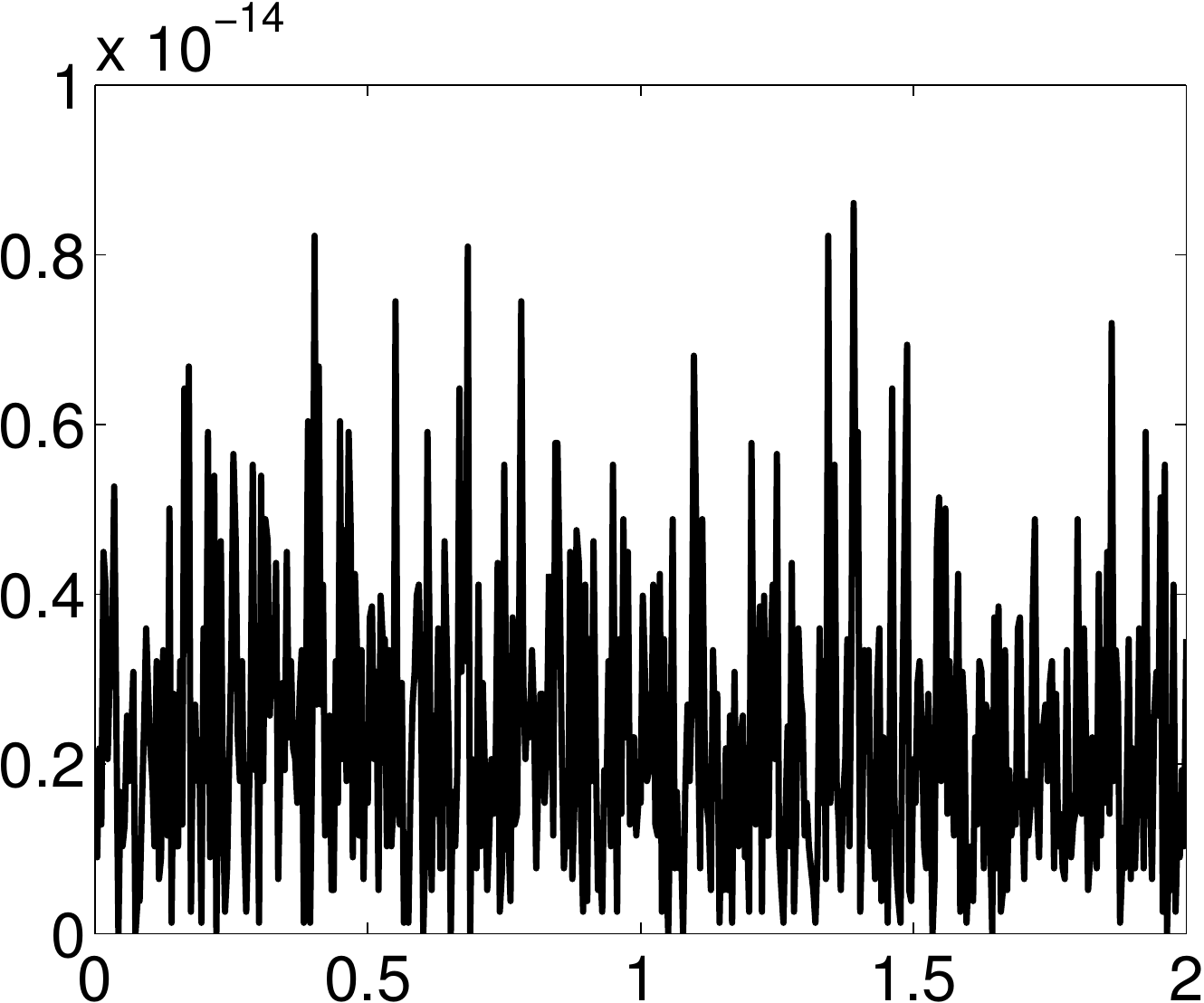}
\caption{Results for Example 4. The relative error in the total surfactant mass versus time for mesh size $h=2/64=0.03125$.  \label{fig:conservationLai}}
\end{center}
\end{figure}

The bulk and surface surfactant concentrations on the moving interface at times $t=0.5, 1, 1.5, 2$ are shown for mesh size $h=2/64=0.03125$ in Fig.~\ref{fig:uBLai} and Fig.~\ref{fig:uSLai}, respectively. The time step size is $k=h/8$. We see in Fig.~\ref{fig:conservationLai} that the proposed method accurately conserves the total surfactant mass. The method in~\cite{ChLai14} extends the bulk equation from $\Omega_1$ to the whole domain by a regularized indicator function, therefore there is a mass leakage to the domain $\Omega_2$ of the order of the regularization parameter. 

\begin{figure}
\begin{center} 
\includegraphics[width=0.45\textwidth]{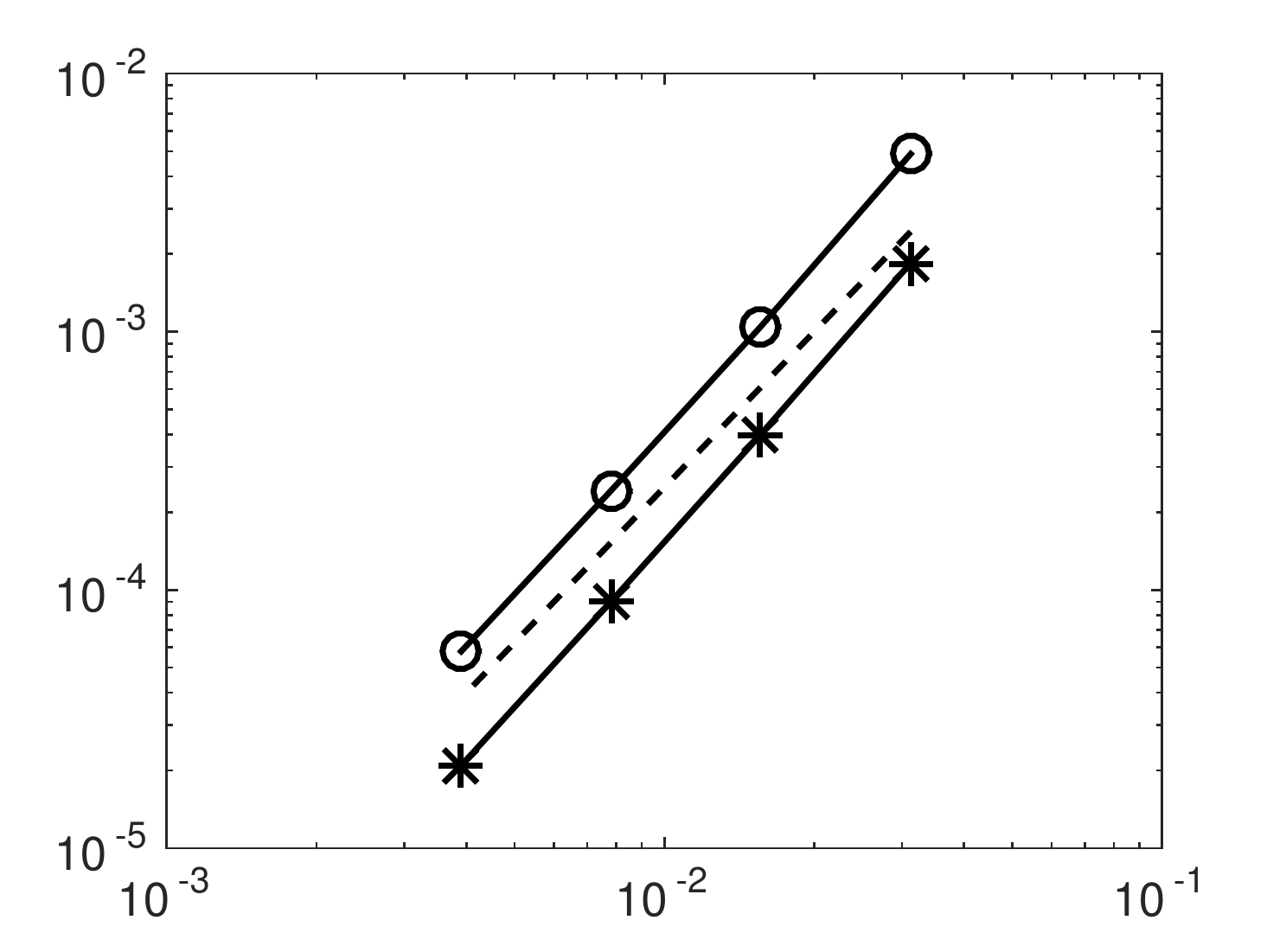}
\caption{Results for Example 4. 
The error $\| (u_{B,h} - u_{B,2h}) \|_{\Omega_{h,1}(0.5)}$ (circles) and $\| (u_{S,h} - u_{S,2h}) \|_{\Gammah(0.5)}$ (stars) measured in the $L^2$ norm versus mesh size $h$. The time step size $k=0.625h$.
The dashed line is proportional to $h^2$. \label{fig:convLai}}
\end{center}
\end{figure}

In Fig.~\ref{fig:convLai} we show the convergence of $\| (u_{B,h} - u_{B,2h}) \|_{\Omega_{h,1}(0.5)}$  (represented by circles) and $\| (u_{S,h} - u_{S,2h}) \|_{\Gammah(0.5)}$ (represented by stars) in the $L^2$ norm. We observe second order convergence which is optimal since we use linear elements in space.
For the two coarsest meshes the error $\| (u_{B,h} - u_{B,2h}) \|_{\Omega_{h,1}(0.5)}$  reported in~\cite{ChLai14} are slightly smaller then the errors we observe however since the method proposed in~\cite{ChLai14} is only first order accurate the errors we observe decrease faster and for the two finest meshes we obtain smaller errors. However, for the mesh sizes shown in the figure the error $\| (u_{S,h} - u_{S,2h}) \|_{\Gammah(0.5)}$ reported in~\cite{ChLai14} is smaller than the error we obtain. This could be explained by the fact that the interface approximation is more accurate in~\cite{ChLai14} with a set of Lagrangian markers. 

\begin{figure}
\begin{center} 
\includegraphics[width=0.4\textwidth]{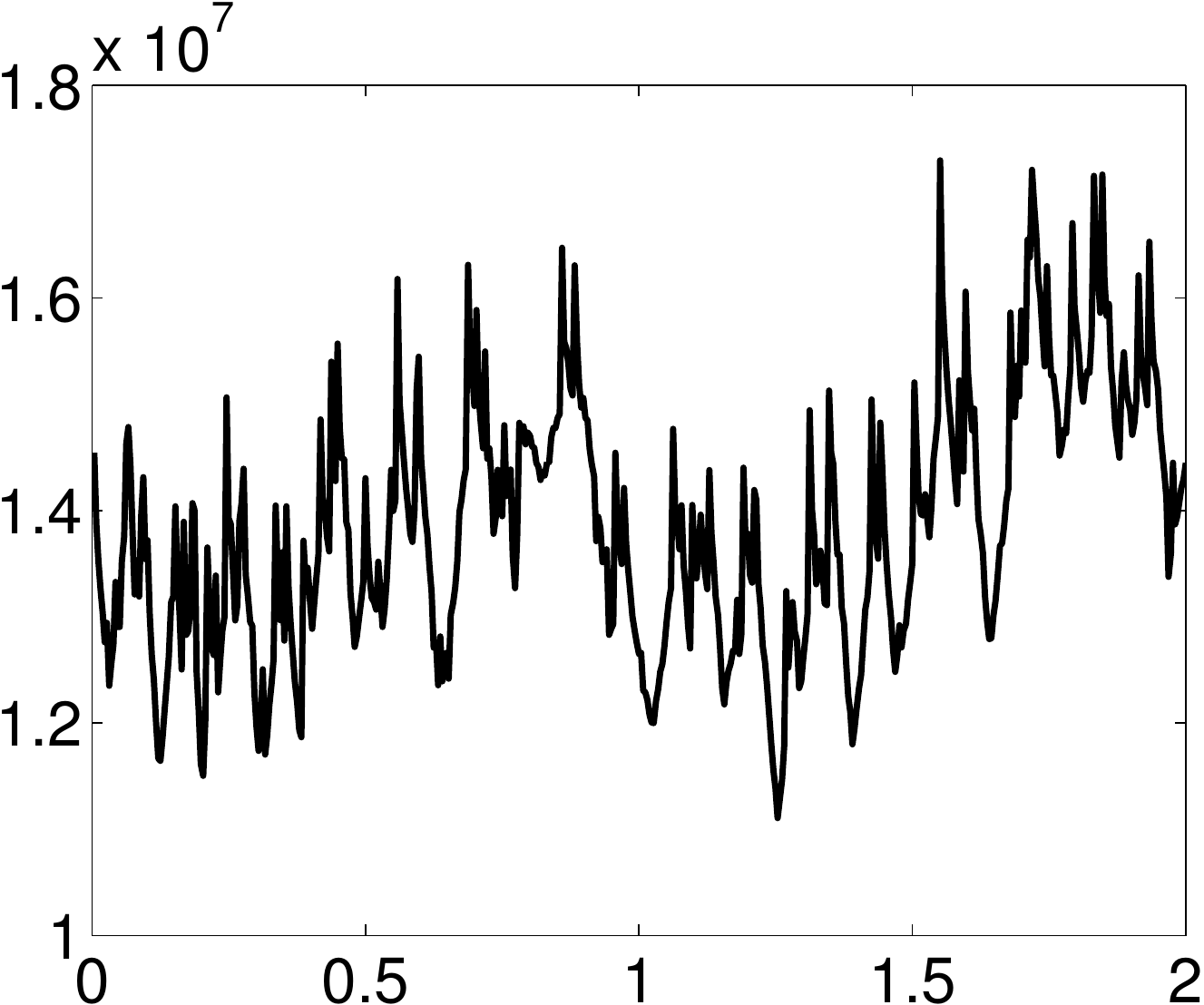}
\caption{Results for Example 4. Condition number versus time for mesh size $h=2/64=0.03125$ and time step size $k=h/8$. 
\label{fig:condvstL}}
\end{center}
\end{figure}

We see in Fig.~\ref{fig:condvstL} where the condition number is shown as a function of time that as the interface evolves the condition number is bounded, independently of how the interface cuts through the mesh.

\section{Conclusion}\label{sec:conc}
We have presented a new finite element method for the coupled bulk-surface problem modeling the evolution of surfactants. The model is given by a time-dependent convection-diffusion equation on the surface (or interface) coupled with a time-dependent convection-diffusion equation in the bulk. The concentration of surfactants in the bulk and the concentration of surfactants on the surface are coupled through a nonlinear term. For this problem we have proposed a space-time CutFEM. The variational formulation is written over the space-time domain. The space-time domain is divided into space-time slabs. Interpolation functions that are continuous in space but discontinuous from one space-time slab to another have been used and the discrete equations have been solved one space-time slab at a time. The space-time formulation allows easily for using unstructured meshes in time slabs which is useful when adaptive schemes are used. For multiphase flow problems adaptive schemes are of great interest because often quantities of interest are on the interface and errors are typically large close to the interface.  

We have added consistent stabilization terms in the weak formulation which guarantee that 1) the system matrices have bounded condition number independently of the interface position relative to the background grid 2) the method is stable in case the problem is convection dominated. A great advantage with the presented method is that it is simple to implement both in two and three space dimensions since it relies on spatial discretizations of the geometry at quadrature points in time and thus much of the computations are the same as in the case of stationary problems. The numerical results demonstrate that the proposed method has optimal convergence order and conserves the total mass of surfactants independently of how the interface cuts through the fixed background mesh.

In this paper, we have used a level set method and represented the interface with linear segments. 
The presented method can be used with any interface representation technique and it is a subject of future work to use a better approximation of the interface with the proposed method.  Also, the velocity field has been given analytically in this paper. However, in future work we will couple the method to a flow solver solving for example the Stokes or the Navier-Stokes equations. 


\section*{Acknowledgements}
This research was supported in part by the Swedish Research Council Grants Nos.\ 2011-4992, 2013-4708, and 2014-4804, the Swedish Foundation for Strategic Research Grant No.\ AM13-0029, and Swedish strategic research programme eSSENCE.

\section{Appendix}\label{sec:appendix}
In some of the numerical examples the non-dimensional form of the equations have been used.  Here we formulate the transport equations in non-dimensional form and report the changes in the weak formulation.   
\subsection{Non-dimensional form of the transport equations for soluble surfactants}
In the following we use the same non-dimensional variables as in~\cite{GaTo12}. Let $L$, $\beta^\infty$, $u_S^\infty$, and $u_B^\infty$ be the characteristic values for length, velocity, and surface and bulk surfactant concentration. The non-dimensional form of the surfactant concentration equations is given by 

\begin{alignat}{2}
\partial_t u_B + \bfbeta \cdot \nabla u_B  - \nabla \cdot \left(\frac{1}{Pe} \nabla u_B \right) &= 0  
\qquad &&\text{in $I\times \Omega_1(t)$}   \label{eq:uBPN} \\ 
-\bfn \cdot \frac{1}{\textrm{Pe}} \nabla u_B &= \textrm{Da}f_{\textrm{coupling}}  
\qquad &&\text{on $\Gamma(t)$} \label{eq:BC1N} \\ 
-\bfn_{\partial \Omega} \cdot \frac{1}{\textrm{Pe}} \nabla u_B &= 0 
\qquad &&\text{on $\partial \Omega$} \label{eq:BC2N} \\
\partial_t u_S + \bfbeta \cdot \nabla u_S +(\divs \bfbeta) u_S  - \divs \left(\frac{1}{\textrm{Pe}_\textrm{S}}\nablas u_S\right) &=   f_{\textrm{coupling}}
\qquad &&\text{on $I\times \Gamma(t)$} \label{eq:uSPN}
\end{alignat}
with
\begin{equation}\label{eq:coupnondim} 
f_{\textrm{coupling}}=\alpha u_B(1-u_S)-\textrm{Bi}u_S
\end{equation} 
Here the non-dimensional numbers Pe and Pe$_{\textrm{s}}$ are the bulk and surface (interfacial) Peclet numbers, defined with respect to the bulk and surface diffusivity $k_B$ and $k_S$, Da is the Damk\"{o}hler number and Bi is the Biot number~\cite{GaTo12} 
\begin{equation}
\textrm{Pe}=\frac{\beta^{\infty}L}{k_B}, \quad \textrm{Pe}_\textrm{S}=\frac{\beta^{\infty}L}{k_S}, \quad \textrm{Da}=\frac{u_S^\infty}{Lu_B^\infty}, \quad  \textrm{Bi}=\frac{k_dL}{\beta^{\infty}}, \quad \alpha=\frac{k_aLu_B^\infty}{\beta^\infty}
\end{equation}

Due to the non-dimensionalization, the conserved total amount of surfactants is 
\begin{equation}\label{eq:conservN}
\int_{\Omega_1(t)} u_B dv +\textrm{Da}\int_{\Gamma(t)} u_S ds
\end{equation}


\subsection{Weak form}
We can write $f_{\textrm{coupling}}$ given in equation~\eqref{eq:coupnondim} in the form 
\begin{equation} 
f_{\textrm{coupling}}={b_B} u_B - {b_S} u_S -b_{BS}u_Bu_S
\end{equation} 
with $b_B=b_{BS}=\alpha$ and $b_S=\textrm{Bi}$. Assuming that ${b_B}$, ${b_S}$, and $\textrm{Da}$ are positive constants we can multiply the bulk PDE~\eqref{eq:uBPN} with a test function $\frac{b_B}{\textrm{Da}}v_B\in H^1(\Omega_1(t))$ and the surface PDE~\eqref{eq:uSPN} as before with a test function $b_Sv_S\in H^1(\Gamma(t))$. This yields after integration by parts and incorporating the boundary conditions the following weak form 
\begin{equation}
 \frac{b_B}{\textrm{Da}}(\partial_t u_B,v_B)_{\Omega_1(t)} + b_S(\partial_t u_S,v_S)_{\Gamma(t)}+a(t,u,v)-(b_{BS}u_Bu_S, b_Bv_B-b_Sv_S)_{\Gamma(t)} =0 \quad \forall  v=(v_B,v_S)\in W
\end{equation}
with
\begin{equation}\label{eq:defaN}
a(t,u,v)= \frac{b_B}{\textrm{Da}}a_B(t,u_B,v_B)+b_Sa_S(t,u_S,v_S)+a_{BS}(t,u,v)
\end{equation}
$a_{BS}$ as before, $a_B$ and $a_S$ as before but with $k_B$ and $k_S$ replaced by Pe and Pe$_{\textrm{s}}$, respectively, see equation~\eqref{eq:contforma}. Note that compared to equation~\eqref{eq:weakform} there is only a change in the coefficient in front of $a_B$ and $(\partial_t u_B,v_B)_{\Omega_1(t)}$ (terms coming from the bulk PDE). 
The same changes are applied to the weak formulation of our space-time CutFEM described in Section~\ref{sec:weakformST}. In addition, the terms containing $\lambda$ and $\mu$ in equation~\eqref{eq:spacetimeform} change since condition~\eqref{eq:conserv} changes to \eqref{eq:conservN}.

\end{document}